\documentclass[11pt]{article}
\usepackage[english]{babel}
\usepackage{times}
\usepackage{draftcopy}
\usepackage{amsmath,amstext,amssymb,amsthm, amsfonts}
\usepackage{graphics,color}
\usepackage{cite}
\newcommand{\bi}{\begin{itemize}}  
\newcommand{\ei}{\end{itemize}}     
\newcommand{\bc}{\begin{center}}  
\newcommand{\ec}{\end{center}}     

\newcommand{\ls}[1]
   {\dimen0=\fontdimen6\the\font \lineskip=#1\dimen0
   \advance\lineskip.5\fontdimen5\the\font \advance\lineskip-\dimen0
   \lineskiplimit=.9\lineskip \baselineskip=\lineskip
   \advance\baselineskip\dimen0 \normallineskip\lineskip
   \normallineskiplimit\lineskiplimit \normalbaselineskip\baselineskip
   \ignorespaces }

\numberwithin{equation}{section}
\newcommand{\slim} {\mathop{\rm lim\,sup}}
\newcommand{\ilim} {\mathop{\rm lim\,inf}}

\newtheorem{lemma}{Lemma}[section]
\newtheorem{theorem}[lemma]{Theorem}
\newtheorem{corollary}[lemma]{Corollary}
\newtheorem{definition}[lemma]{Definition}
\newtheorem{proposition}[lemma]{Proposition}
\newtheorem{example}[lemma]{Example}

\newtheorem{remark}[lemma]{Remark}
\newtheorem{Assumption}[lemma]{Assumption}

\def\S{\mathbb{S}}

\def\K{\mathbb{K}}

\def\Y{\mathbb{Y}}

\def\X{\mathbb{X}}

\def\A{\mathbb{A}}

\def\R{\overline{\mathbb{R}}}
\def\P{\mathbb{P}}

\def\B{\mathcal{B}}

\def\bb{\mathbb{B}}

\def\lv(#1){{\underline{#1}}}
\def\uv(#1){{\overline{#1}}}

\def\lv(#1){{\underline{#1}}}
\def\uv(#1){{\overline{#1}}}
\def\uh(#1){{\hat{#1}}}
\def\wsol(#1){{\rm Sol}_w(#1)}
\def\sol(#1){{\rm Sol}(#1)}
\def\esol(#1){{\rm Sol}_\varepsilon(#1)}

\title{ Solutions for Zero-Sum Two-Player Games with Noncompact Decision Sets and Unbounded Payoffs}

\begin{document}

\maketitle

\begin{center}
Eugene~A.~Feinberg \footnote{Department of Applied Mathematics and
Statistics,
 Stony Brook University,
Stony Brook, NY 11794-3600, USA, eugene.feinberg@sunysb.edu},\
Pavlo~O.~Kasyanov  \footnote{Institute for Applied System Analysis,
National Technical University of Ukraine ``Igor Sikorsky Kyiv Polytechnic
Institute'', Peremogy ave., 37a, build, 35, 03056, Kyiv, Ukraine,\
kasyanov@i.ua.},\
and Michael~Z.~Zgurovsky \footnote{National Technical University of Ukraine
``Igor Sikorsky Kyiv Polytechnic
Institute'', Peremogy ave., 37, build, 1, 03056, Kyiv, Ukraine,\
mzz@kpi.ua.}

\bigskip
\end{center}

\begin{abstract}
This paper provides sufficient conditions for the existence of solutions for two-person zero-sum games with  inf/sup-compact payoff functions and with possibly noncompact decision sets for both players. Payoff functions may be unbounded, and we do not assume any convexity/concavity-type conditions.  For such games expected payoff may not exist for some pairs of strategies.  The results of this paper imply several classic {facts}.  The paper also provides sufficient conditions for the existence of a value and solutions for each player.   The results of this paper are illustrated with the number guessing game.
\end{abstract}


\section{Introduction}\label{intro}

This paper studies two-person zero-sum games with possibly noncompact decision sets and unbounded payoff functions, which may not satisfy any convexity or concavity conditions. The payoff functions may be bounded neither from below nor from above. The players can use mixed strategies.

Since the payoff function is unbounded from above and from below,  the expected payoffs are uncertain for some pairs of mixed strategies. This paper provides sufficient conditions for the existence of a solution or, in other words, for the existence of a pair of equilibrium strategies; see Theorem~\ref{teor:mainonestep} below.   Each strategy in an equilibrium pair satisfies the following property: if the corresponding player chooses this strategy, then the expected finite or infinite payoff is defined for every strategy played by the opponent.  Strategies of each player with this property are called { sensible}.  Of course, if the payoff functions are bounded from below or from above, then every strategy is { sensible}.  However, if the decision sets for both players are noncompact, then natural mild continuity conditions (see Assumptions~\ref{AsMain}(a2,b2) below) imply that the payoff function is unbounded from below and from above.   This paper also provides sufficient conditions for the existence of values and solutions for each player.

The currently available general result on the existence of solutions for infinite games, Mertens et. al.~\cite[Theorem I.2.4]{MSZ}, assumes that  both decision sets are compact and the payoff function is bounded from below or above.  Its proof { in \cite{MSZ}} uses Sion's minimax theorem  for quasiconvex/quasiconcave functions applied to the sets of mixed strategies of the original game. This method  works only if the payoff function is bounded from below or from above.  If the payoff function is bounded neither from below nor from above, the expected payoffs are not defined for some pairs of mixed strategies.  Sion's minimax theorem assumes that  decision sets are subsets of linear topological spaces, and one of them is compact. There are minimax theorems for more general spaces, than linear topological spaces  \cite{Khanh, Park}, but they also assume that  one of decision sets is compact; see, e.g., \cite[Corollary 4.1]{Khanh}.

Our proofs are based on the lopsided minimax equality established in \cite[Theorem 18]{OSCG} for games with two possibly noncompact decision sets; see
Theorem~\ref{th:exvalPaper1}(a) below. The classic sufficient conditions for the validity of such minimax equalities require that one of the decision sets is compact; see Mertens et al.~\cite[Propositions I.1.9 and I.2.2]{MSZ} and Alpern and Gal~\cite{AG}.  Theorem~\ref{th:exvalPaper1} leads to Assumptions~{\ref{ass:MAIN-NEW-A} and \ref{ass:MAIN-NEW-B}} introduced in Section{s~\ref{sec:2} and \ref{Sec:MR}}.  These assumptions  are equivalent to the assumption that, if one of the players chooses an { unsensible} strategy, another player can respond with a strategy causing the infinite loss to the player who chose the { unsensible} strategy. The lopsided minimax equality  in \cite[Theorem 18]{OSCG} (see Theorem~\ref{th:exvalPaper1}(a) below), is used in \cite[Theorem 20]{OSCG} to prove the existence of a solution if payoffs are bounded from below.

Section~\ref{sec:2} of this paper provides main definitions and  preliminary results, some of which are taken from our previous paper \cite{OSCG}.  Section~\ref{Sec:MR} provides main results.  In addition to Assumptions~{\ref{ass:MAIN-NEW-A} and \ref{ass:MAIN-NEW-B}} mentioned above, the existence of a solution  is established in Theorem~\ref{teor:mainonestep} under the conditions that the payoff function is lower/upper semi-continuous  in the corresponding variables, and  the payoff function is inf/sup-compact in the decision variable corresponding to decisions of the player for at least one decision chosen by another player.  The same conditions are assumed for the payoff function  in Aubin and Ekeland~\cite[Theorem 6.2.8]{ObEk} stating the existence of solutions within the class of pure strategies for convex/concave payoff functions.  The results are illustrated in Section~\ref{sec:ex} in which the following game is considered.  Two players choose nonnegative numbers, and the payoff is a polynomial function of the difference between these numbers.  We completely classify this game:  in some cases there are no solutions, in some cases solutions exist, an in some cases solutions exist, but pure solutions do not exist. Some of the proofs are provided in Section~\ref{sec:proofs}.

Our initial motivation for studying games with unbounded payoffs in Feinberg et al.~\cite{OSCG} and in this paper was originated by the progress in the theory of Markov Decision Processes with possibly noncompact decision sets and unbounded costs
that led to the extension of Berge's maximum theorem \cite[p.116]{Ber} to possibly noncompact decision sets; see Feinberg et al.~\cite{FKK, FKV, Feinberg_et_al}, Feinberg and Kasyanov~\cite{FKSVAN}.  These results were applied in \cite{OSCG} to games with perfect information and, as described above, some results for games with simultaneous moves are also obtained in \cite{OSCG}. In particular, the existence of solutions is established in \cite[Theorem 20]{OSCG} for payoffs bounded below, and this and other assumptions imply compactness of one of the decision sets.
This paper studies more general models, when  both decision sets may be noncompact. These models have potential applications to stochastic games; see  Ja\'shkewicz and Nowak~\cite{Jan, JanS}, Mertens et al.~\cite[Chapter VII]{MSZ}, and references therein for the literature on stochastic games.

\section{Definitions and Preliminary Results}\label{sec:2}

 Let  $\S$ be a metric space, and ${\mathcal B}(\S)$ be the Borel
$\sigma$-field on $\S,$ that is, the $\sigma$-field is generated by all
open sets of the metric space $\S.$ For a nonempty Borel subset $S\subset \S,$ let
${\mathcal B}(S)$ denote the $\sigma$-field whose elements are
intersections of $S$ with elements of ${\mathcal B}(\S).$  Observe that $S$ is a metric space with the same metric as on $\S.$  Therefore,
${\mathcal B}(S)$ is its Borel $\sigma$-field. Let $\P(\S)$ be the set of probability measures on
$(\S,{\mathcal B}(\S)).$ We denote by $\P^{fs}(S)$  the set of all probability measures whose supports are finite subsets of  the set $S.$ A sequence of probability measures
$\{\mu^{(n)}\}_{n=1,2,\ldots}$ from $\P(\S)$ \textit{converges weakly} to $\mu\in\P(\S)$ if for
each bounded continuous function $f$ on $\S$
\[\int_\S f(s)\mu^{(n)}(ds)\to \int_\S f(s)\mu(ds) \qquad {\rm as \quad
}n\to\infty.
\]
We endow $\P(\S)$ with the topology of the weak convergence of probability measures on $\S$. If
$\S$ is
a separable metric space, then $\P(\S)$ is separable metric space too and the set $\P^{fs}(\S)$ is dense in $\P(\S)$; Parthasarathy \cite[Chapter II, Theorems 6.2 and 6.3]{Part}. Let $\R:=\mathbb{R}\cup\{\pm\infty\},$  where $\mathbb{R}$ is the set of real numbers.

 An integral $\int_\S f (s)\mu(ds)$ of a measurable ${\R}$-valued function $f$ on $\S$ over the measure $\mu \in
\P(\S)$ is well-defined if either $-\infty<\int_\S f^- (s)\mu(ds)$ or $\int_\S f^+ (s)\mu(ds)<+\infty, $  where $f^-(s)=\min\{f(s),0\},$ $f^+(s)=\max\{f(s),0\},$ $s\in \S.$  If the integral is well-defined, then $\int_\S f (s)\mu(ds):=\int_\S f^+ (s)\mu(ds)+ \int_\S f^-(s)\mu(ds).$

\begin{definition}\label{defi:game}  A \textit{  two-person zero-sum game} is a triplet $\{\A,\bb, c\},$ where
  \begin{itemize}
\item[{\rm(i)}] $\A$ is the \textit{space of decisions for Player I}, which is a nonempty Borel subset of a Polish space;
\item[{\rm(ii)}] $\bb$ is the \textit{space of decisions for Player II}, which is a nonempty Borel subset of a Polish space;
\item[{\rm(iii)}] the \textit{payoff} from Player~I to  Player II, $-\infty< c(a,b)< +\infty,$ for choosing decisions $a\in \A$ and $b\in \bb,$ is a \textit{measurable} function on $\A\times\bb;$
    \item[{\rm(iv)}] for each $b\in\bb$ the function $a\mapsto c(a,b)$ is bounded from below on $\A;$
\item[{\rm(v)}]  for each $a\in\A$ the function $b\mapsto c(a,b)$ is bounded from above on $\bb.$
\end{itemize}
\end{definition}
\textit{The game is played as follows}:

$\bullet$ a decision-makers (Players I and II) choose simultaneously respective decisions $a\in \A$ and $b\in\bb;$

$\bullet$ the result $(a,b)$ is announced to both of them;

$\bullet$ Player I pays Player II the amount $c(a,b).$

 Everywhere in this paper, except in some pathological situations described in Section~\ref{sec:ex}, we assume that a game $\{\A,\bb, c\}$ satisfies conditions (i-v) from Definition~\ref{defi:game} and consider only such games.

\textit{Strategies} (sometimes called ``mixed strategies'') for Players I and II are probability measures 
$\pi^\A\in \P( \mathbb{A})$ and $\pi^\bb\in \P( \mathbb{B})$ respectively.
 Moreover, $\pi^\A$ ($\pi^\bb$)
is called \textit{pure} if the probability measure
$\pi^\A(\,\cdot\,)$ ($\pi^\bb(\,\cdot\,)$) is concentrated at one point.
Note that $\P(\A)$
is the
\textit{set of   strategies} for Player I,
and $\P(\bb)$ is the
\textit{set of   strategies} for Player II.

With a slight abuse of notation, we identify a pure strategy with the decision it chooses.  In particular, $\A$ and $\bb$ are the sets of pure strategies for Players I and II respectively.  We  usually write $a$ instead of $\delta_{\{a\}}$ and $b$ instead of  $\delta_{\{b\}},$  where  $\delta_{\{a\}}$ and $\delta_{\{b\}}$ are probability measures on $(\A,{\mathcal B}(\A))$ and $(\bb,{\mathcal B}(\bb))$ concentrated at the points $a\in\A$ and $b\in\bb$ respectively.

Let 
\begin{equation*}
\uh(c)^\oplus(\pi^\A,\pi^\bb):=\int_{\A}\int_{\bb}c^+(a,b)\pi^\bb(db)\pi^\A(da),\quad
\uh(c)^\ominus(\pi^\A,\pi^\bb):=\int_{\A}\int_{\bb}c^-(a,b)\pi^\bb(db)\pi^\A(da)
\end{equation*}
for each $(\pi^\A,\pi^\bb)\in \P(\A)\times\P(\bb).$ Then
the \textit{expected
payoff} from Player~I to  Player~II is 
\begin{equation}\label{EFnC}
\uh(c)(\pi^\A,\pi^\bb):= \uh(c)^\oplus(\pi^\A,\pi^\bb)+\uh(c)^\ominus(\pi^\A,\pi^\bb),
\end{equation}
and it is well-defined, if either $\uh(c)^\oplus(\pi^\A,\pi^\bb)<+\infty$ or $-\infty<\uh(c)^\ominus(\pi^\A,\pi^\bb),$ where $(\pi^\A,\pi^\bb)\in \P(\A)\times\P(\bb).$  Of course, 
when the function $c$ is unbounded from below and from above, the quantity $\uh(c)(\pi^\A,\pi^\bb)$ is  \textit{undefined} for some $(\pi^\A,\pi^\bb)\in\P(\A)\times\P(\bb).$
Assumptions~(iv) and (v) for the game $\{\A,\bb,c\}$ imply that
 $-\infty<\uh(c)^\ominus(\pi^\A,b)$ for each $\pi^\A\in \P(\A)$ and $b\in\bb,$ and
 $\uh(c)^\oplus(a,\pi^\bb)<+\infty$ for each $a\in\A$ and $\pi^\bb\in \P(\bb)$ respectively.

 The set of   strategies for each player is partitioned into the sets of \textit{{ sensible} strategies} $\P^{S}(\A)$ and $\P^{S}(\bb)$  (strategies, for which the expected payoff is well-defined for all strategies played by another player) and \textit{{ unsensible} strategies} $\P^{U}(\A)$ and $\P^{U}(\bb)$:
 \[
\begin{aligned}
\P^{S}(\A)&:=\{\pi^\A\in\P(\A)\,:\, \uh(c)(\pi^\A,\pi^\bb_*)\mbox{ is well-defined for each }\pi^\bb_* \in \P(\bb)\}, \\
&\P^{U}(\A):=\{\pi^\A\in\P(\A)\,:\, \uh(c)(\pi^\A,\pi^\bb_*)\mbox{ is undefined for some }\pi^\bb_* \in \P(\bb)\},\\
\P^{S}(\bb)&:=\{\pi^\bb\in\P(\bb)\,:\, \uh(c)(\pi^\A_*,\pi^\bb)\mbox{ is well-defined for each }\pi_*^\A\in\P(\A)\},\\
&\P^{U}(\bb):=\{\pi^\bb\in\P(\bb)\,:\, \uh(c)(\pi^\A_*,\pi^\bb)\mbox{ is undefined for some }\pi_*^\A\in\P(\A)\}.
\end{aligned}
\]
Assumptions~(iv) and (v) from Definition~\ref{defi:game} imply that $\P^{fs}(\A)\subset \P^{S}(\A)$ and $\P^{fs}(\bb)\subset \P^{S}(\bb)$. Therefore, $\P^{S}(\A)$ is dense in $\P(\A)$ and
$\P^{S}(\bb)$ is dense in $\P(\bb).$

Let us introduce  the \textit{worst gain} of Player~II for choosing a strategy $\pi^\bb\in \P(\bb)$ and \textit{worst loss} of Player~I for choosing a strategy $\pi^\A\in\P(\A)$ respectively:
\begin{equation}\label{eq:maxmin}
\uh(c)^{\flat}(\pi^\bb):=\inf\limits_{a\in \A } \uh(c)(a,\pi^\bb)\quad\mbox{and}\quad \uh(c)^{\sharp}(\pi^\A):=\sup\limits_{b\in \bb} \uh(c)(\pi^\A,b).
\end{equation}
 These definitions and  assumptions (iv,v) from Definition~\ref{defi:game} of the game $\{\A,\bb,c\}$ imply
\begin{equation}\label{eq:safe00}
-\infty<\uh(c)^{\flat}(b)\le c(a,b)\le \uh(c)^{\sharp}(a)<+\infty,\qquad\qquad a\in\A,\ b\in\bb,
\end{equation}
According to Feinberg et al. \cite[Theorem~17 and Lemma 6]{OSCG}, 
\begin{equation}\label{eq:prop:mainIaIb}
\uh(c)^{\flat}(\pi^\bb)=\inf_{\pi_*^\A\in \P^S(\A)} \uh(c)(\pi_*^\A,\pi^\bb)\quad\mbox{and}\quad
\uh(c)^{\sharp}(\pi^\A)=\sup_{\pi_*^\bb\in\P^S(\bb)} \uh(c)(\pi^\A,\pi_*^\bb ),
\end{equation}
the  inequality
\begin{equation}\label{eq:g}
\uh(c)^{\flat}(\pi^\bb)\le \uh(c)^{\sharp}(\pi^\A)
\end{equation}
holds for each $\pi^\A\in\P(\A)$ and $\pi^\bb\in\P(\bb)$ such that $\uh(c)(\pi^\A,\pi^\bb)$ is well-defined,
the functions $\uh(c)^{\sharp}$ and $\uh(c)^{\flat}$ are convex on $\P(\A)$
and concave on $\P(\bb)$ respectively, and the following level sets
\[
\begin{aligned}
\P^\sharp_\alpha(\A):=\{\pi_*^\A\in \P(\A)\,:\, \uh(c)^\sharp(\pi^\A_*)\le \alpha\},\quad \P^\flat_\alpha(\bb):=\{\pi_*^\bb\in \P(\bb)\,:\, \uh(c)^\flat(\pi_*^\bb)\ge \alpha\},
\end{aligned}
\]
are convex for all $\alpha\in\R.$

In particular, $\uh(c)(\pi^\A,\pi^\bb)$  is well-defined for $\pi^\A\in\P^S(\A)$ and $\pi^\bb\in\P^S(\bb).$  Therefore, inequality \eqref{eq:g} holds for $\pi^\A\in\P^S(\A)$ and $\pi^\bb\in\P^S(\bb).$ It is not clear whether inequality \eqref{eq:g} holds for $\pi^\A\in\P(\A)$ and $\pi^\bb\in\P(\bb)$  when $\uh(c)(\pi^\A,\pi^\bb)$ is undefined.

 The following definition
 introduces the lower and upper values  in a slightly more delicate way than it is usually done when $\uh(c)(\pi^\A,\pi^\bb)$ are always defined.

\begin{definition}\label{Defulvalgen}{\rm (Lower and upper values  of the game).
The quantities
\begin{equation}\label{eq:lowupvalB}
v ^\flat:=\sup_{\pi^\bb\in\P^S(\bb)}\uh(c)^{\flat}(\pi^\bb)\quad\mbox{and}\quad v ^\sharp:=\inf_{\pi^\A\in\P^S(\A)} \uh(c)^{\sharp}(\pi^\A)
\end{equation}
are the \textit{lower and upper values of the game} $\{\A,\bb,c\}$ respectively.}
\end{definition}
Since inequality \eqref{eq:g} is valid for every pair of strategies $\pi^\A\in\P^S(\A)$ and $\pi^\bb\in\P^S(\bb),$ the inequality $v^\flat \le v^\sharp$ holds.  In addition,
$-\infty< \sup_{b\in\bb}\uh(c)^{\flat}(b)\le \sup_{\pi^\bb\in\P^S(\bb)}\uh(c)^{\flat}(\pi^\bb)=
v^\flat,$ where the second inequality holds because every pure strategy is { sensible} in view of assumptions (iv) and (v) from   Definition~\ref{defi:game}, and the first inequality follows from \eqref{eq:safe00}.  Similarly, $v ^\sharp<+\infty.$   Therefore,
\begin{equation}\label{eq2.7aEF}
-\infty<v^\flat \le v^\sharp<+\infty.
\end{equation}
Equalities (\ref{eq:prop:mainIaIb}) and   (\ref{eq:lowupvalB}) imply
\begin{equation}\label{eq:newstar}
v ^\flat=\sup_{\pi^\bb\in\P^S(\bb)}\inf_{\pi^\A\in\P^S(\A)}\uh(c)(\pi^\A,\pi^\bb)\quad\mbox{and}\quad v ^\sharp=\inf_{\pi^\A\in\P^S(\A)}\sup_{\pi^\bb\in\P^S(\bb)} \uh(c)(\pi^\A,\pi^\bb).
\end{equation}
If $\uh(c)(\pi^\A,\pi^\bb)$ is well-defined for each $(\pi^\A,\pi^\bb)\in\P(\A)\times \P(\bb),$ as this takes place when $c$ is bounded from below or above on $\A\times\bb,$ then  the upper and lower values of the game defined in (\ref{eq:lowupvalB})
coincide with their classic definition, that is,
\begin{equation}\label{eq:lowupvalBCL}
v ^\flat=\sup_{\pi^\bb\in\P(\bb)}\inf_{\pi^\A\in\P(\A)}\uh(c)(\pi^\A,\pi^\bb)\quad\mbox{and}\quad v ^\sharp=\inf_{\pi^\A\in\P(\A)}\sup_{\pi^\bb\in\P(\bb)} \uh(c)(\pi^\A,\pi^\bb).
\end{equation}

%
\begin{definition}\label{defi:soleachequil}{\rm (Solution for a player).
 A strategy $\pi^\A\in \P^S(\A)$  ($\pi^\bb\in \P^S(\bb)$) is called \emph{a  solution for Player~ I (II)} if
\begin{equation}\label{eq:solutionI or II}
\uh(c)^\sharp(\pi^\A)=v^\sharp \quad (\uh(c)^\flat(\pi^\bb)=v^\flat).
\end{equation}
A solution $\pi^\A$  ($\pi^\bb$) for Player I (II)} is  called \textit{pure} if the respective strategy
$\pi^\A$ ($\pi^\bb$) is pure.
\end{definition}
\begin{definition}\label{Def5.8}{\rm (Value of a game).}
If the equality
\begin{equation}\label{eq:valval}
v^\flat =v^\sharp
\end{equation}
holds, then we say that this common quantity is the \textit{value
of the game $\{\A,\bb,c\}$}.
We denote the value by $v$.
\end{definition}
In view of this definition, if the value $v$ exists, then it is unique and, in view of \eqref{eq2.7aEF}, $v$ is a real number. If $\uh(c)(\pi^\A,\pi^\bb)$ is well-defined for each $(\pi^\A,\pi^\bb)\in\P(\A)\times \P(\bb),$ as this takes place if and only if $c$ is bounded from below or above on $\A\times\bb,$ then  the value defined in (\ref{eq:valval})
coincides with its classic definition, that is,
\[
\begin{aligned}
v=v^\flat =\sup_{\pi^\bb\in\P(\bb)}\uh(c)^{\flat}(\pi^\bb)=\sup_{\pi^\bb\in\P(\bb)}\inf_{\pi^\A\in\P(\A)}\uh(c)(\pi^\A,\pi^\bb)=
v^\sharp =\inf_{\pi^\A\in\P(\A)} \uh(c)^{\sharp}(\pi^\A)=\inf_{\pi^\A\in\P(\A)}\sup_{\pi^\bb\in\P(\bb)}\uh(c)(\pi^\A,\pi^\bb).
\end{aligned}
\]
\begin{lemma}\label{lemm:valdefi}
If the value $v$ of  a game $\{\A,\bb,c\}$ exists, then 
\begin{equation}\label{eq:valvalclass}
v=\sup_{\pi^\bb\in\P(\bb)}\uh(c)^{\flat}(\pi^\bb)=\inf_{\pi^\A\in\P(\A)} \uh(c)^{\sharp}(\pi^\A).
\end{equation}
\end{lemma}
\begin{proof}
Note that
\begin{equation}\label{eq:ineqnew1}
v^\flat=\sup_{\pi^\bb\in\P^S(\bb)}\uh(c)^{\flat}(\pi^\bb)\le \sup_{\pi^\bb\in\P(\bb)}\uh(c)^{\flat}(\pi^\bb)\le\inf_{\pi^\A\in\P^S(\A)} \uh(c)^{\sharp}(\pi^\A)=v^\sharp,
\end{equation}
where the equalities follow from the definitions of lower and upper values (\ref{eq:lowupvalB}), the first inequality follows from $\P^S(\bb)\subset \P(\bb),$ and the second inequality follows from inequality (\ref{eq:g})
because $\uh(c)(\pi^\A,\pi^\bb)$ is well-defined for each pair $(\pi^\A,\pi^\bb)\in\P^S(\A)\times\P(\bb).$ Similarly,
\begin{equation}\label{eq:ineqnew2}
v^\flat=\sup_{\pi^\bb\in\P^S(\bb)}\uh(c)^{\flat}(\pi^\bb)\le \inf_{\pi^\A\in\P(\A)} \uh(c)^{\sharp}(\pi^\A)\le \inf_{\pi^\A\in\P^S(\A)} \uh(c)^{\sharp}(\pi^\A)=v^\sharp.
\end{equation}
Therefore, (\ref{eq:valvalclass}) follows from (\ref{eq:ineqnew1}) and (\ref{eq:ineqnew2}) because  $v^\flat=v^\sharp=v.$
\end{proof}

We recall that for a metric space $\S$ a function $f:\S\to\R$ is
called \textit{lower semi-continuous at $s\in \S$}, if for each sequence
$\{s^{(n)}\}_{n=1,2,\ldots},$ that converges to $s$ in $\S,$ the
inequality $\ilim_{n\to\infty} f(s^{(n)})\ge f(s)$ holds. A function
$f:\S\to\R$ is called \textit{upper semi-continuous at
$s\in \S$}, if $-f$ is lower semi-continuous at $s\in \S.$ A function
$f: \S\to\R$ is called \textit{lower / upper semi-continuous}, if $f$ is lower / upper semi-continuous at each $s\in \S.$
A function
$f: \S\to\R$ is called \textit{inf-compact on
$\S$}, if all the level sets $\{\{s\in \S \, : \,  f(s)\le \lambda\}\}_{\lambda\in\mathbb{R}}$ are compact in $\S.$ A function
$f:\S\to\R$ is called \textit{sup-compact on
$\S$} if $-f$ is inf-compact on~$\S.$

The following assumptions contains all the semi-continuity assumptions on the payoff function $c$ used in this paper.
\begin{Assumption}\label{AsMain} The payoff function $c:\X\times\A\to\mathbb R$ satisfies the following conditions:
\begin{itemize}
\item[{\rm(a1)}] for each $b\in \bb$ the function $a\mapsto c(a,b)$ is lower semi-continuous,
\item[{\rm(a2)}] there exists $b_0\in \bb$ such that the function $a\mapsto c(a,b_0)$ is inf-compact on $\A;$
    \item[{\rm(b1)}] for each $a\in \A$ the function $b\mapsto c(a,b)$ is upper semi-continuous;
\item[{\rm(b2)}] there exists $a_0\in \A$ such that the function $b\mapsto c(a_0,b)$ is sup-compact on $\bb.$
\end{itemize}
\end{Assumption}

The following theorem, whose part (a) is taken from Feinberg et al.~\cite[Subsection~5.4]{OSCG}, provides sufficient conditions for the
existence of a value of the game $\{\A,\bb, c \}$ in the sense of Definition~\ref{Def5.8}.

\begin{theorem}\label{th:exvalPaper1} (a) {\rm (Existence of a lopsided value; Feinberg et al.~\cite[Theorem 18 and Corollary 4]{OSCG}).}
 If a  two-person zero-sum game $\{\A,\bb, c \}$   satisfies Assumptions~\ref{AsMain}(a1,a2),
then
\begin{equation}\label{eq:valvalPaper1}
\min_{\pi^\A\in\mathbb{P}(\A)} {\hat{c}}^{\sharp}(\pi^\A)=  \sup_{\pi^\bb\in\mathbb{P}^{fs}(\bb)}{\hat{c}}^{\flat}(\pi^\bb)=v^\flat,
\end{equation}
and the set $\mathbb{P}^\sharp_{v^\flat}(\A)$ is a
nonempty convex compact subset of $\mathbb{P}(\A).$

(b) {\rm (Existence of the value).} Under the assumptions of statement (a), the value  $v=v^\flat=v^\sharp$ of the game exists 
if and only if
\begin{equation}\label{eq:valUvalS}
  v^\sharp=
\inf_{\pi^\A\in\P(\A)}\uh(c)^{\sharp}(\pi^\A).
\end{equation}
\end{theorem}
 \begin{proof}[Proof of Theorem~\ref{th:exvalPaper1}(b)] Let the assumptions of Theorem~\ref{th:exvalPaper1}(a) hold.
In view of \eqref{eq:valvalPaper1}, the infimum in \eqref{eq:valUvalS} can be replaced with the minimum.  If $v^\flat=v^\sharp,$ then \eqref{eq:valvalPaper1} implies \eqref{eq:valUvalS}.  If \eqref{eq:valUvalS} holds, then \eqref{eq:valvalPaper1} implies  $v^\flat=v^\sharp.$
\end{proof}

Theorem~\ref{th:exvalPaper1}(b)  is useful for  proving the existence of the value. Observe that \eqref{eq:valUvalS} holds { under the following condition:
\begin{Assumption}\label{ass:MAIN-NEW-A}
$\uh(c)^\sharp(\pi^\A)> v^\sharp$ for all $\pi^\A\in\P^U(\A).$
\end{Assumption}}
\noindent This is true because $\P(\A)=\P^S(\A)\cup\P^U(\A),$ and, if {Assumption~\ref{ass:MAIN-NEW-A}} is correct, then \eqref{eq:valUvalS} becomes the definition of $v^\sharp$ given in \eqref{eq:lowupvalB}.
Therefore, {Assumption~\ref{ass:MAIN-NEW-A}} and Assumptions~\ref{AsMain}(a1,a2)  imply the existence of the value of the game. {Moreover, if the value $v$ exists and Assumption~\ref{ass:MAIN-NEW-A} holds, then
\begin{equation}\label{eq:PvissafeEF}
\P^\sharp_{v}(\A)\subset \P^S(\A).
\end{equation}}
The similar observations are applicable to the lower value $v^\flat,$ when condition (\ref{eq:safe1a}) is replaced     symmetrically with condition \eqref{eq:safe1b} for Player II  presented below.

{Note that Assumption~\ref{ass:MAIN-NEW-A} holds if
\begin{equation}\label{eq:safe1a}
\uh(c)^\sharp(\pi^\A)=+\infty\mbox{ for all }\pi^\A\in\P^U(\A).
\end{equation}}
However, as Example 3 in \cite{OSCG} demonstrates, it is possible that $\uh(c)^\sharp(\pi^\A)< +\infty$ for some $\pi^\A\in\P^U(\A)$ even under stronger conditions than Assumptions~\ref{AsMain}(a1,a2). We also note that statement (\ref{eq:safe1a}) holds if and only if $\P^\sharp_{\alpha}(\A)\subset \P^S(\A)$ for each $\alpha\in\mathbb{R}.$

\begin{remark}\label{rem:cond(i,ii)}
{\rm Assumptions~\ref{AsMain}(a1,a2) are natural for the existence of a solution for Player~I.  For example, they are assumed in  \cite[Theorem 6.2.7]{ObEk}, where the existence of a solution for Player~I is stated for   payoff functions $c(a,b)$ being convex in $a$ and concave in $b.$  In particular, if the decision set of Player II is a singleton, that is, $\bb=\{b_0\},$ then the game becomes an optimization problem. For this optimization problem, inf-compactness Assumption~\ref{AsMain}(a2) is a natural sufficient condition for the existence of a minimum, and this minimum corresponds to the solution for Player~I and for the game.   If the function $c$ is bounded above, then Assumption~\ref{AsMain}(a2) implies that the set $\A$ is compact because $\A=\{a\in \A: c(a,b_0)\le \lambda\}$ for some $\lambda\in\mathbb{R}.$  If the set $\A$ is compact, then Assumption~\ref{AsMain}(a1) implies Assumption~\ref{AsMain}(a2).}
\end{remark}

Let us consider some corollaries from Theorem~\ref{th:exvalPaper1}.   We recall that, as discussed in Section~\ref{sec:2}, the action sets $\A$ and $\bb$ can be  identified with the sets of pure strategies.
\begin{corollary}\label{cor:DeltaEF} (a) For  is a two-person game  $\{\A,\bb,c\}$,
\begin{equation}\label{eq:sup01052019EF}
\uh(c)^\sharp(\pi^\A)=\sup_{\pi^\bb\in\Delta^\prime(\bb)} \uh(c)(\pi^\A,\pi^\bb),\qquad\qquad \pi^\A\in\P(\A),
\end{equation}
for all $\Delta^\prime(\bb)\subset \P(\bb)$ such that  $\bb\subset \Delta^\prime(\bb)\subset \P^S(\bb).$ 

(b) {\rm (\cite[Corollary 3]{OSCG})}. If Assumptions~\ref{AsMain}(a1,a2) hold, then
\[
\sup_{\pi^\bb\in\Delta(\bb)} \uh(c)^\flat(\pi^\bb)=\min_{\pi^\A\in\P(\A)} \uh(c)^\sharp(\pi^\A),
\]
for all $\Delta(\bb)\subset \P(\bb)$ such that $\P^{fs}(\bb)\subset \Delta(\bb)\subset \P^S(\bb).$
\end{corollary}
\begin{proof}
Statement (a) follows from \eqref{eq:maxmin} and \eqref{eq:prop:mainIaIb}.
\end{proof}
The following corollary states the classic minimax equality, which  is well-known for a compact set $\A;$ see {\rm Mertens et al.~\cite[Proposition I.1.9]{MSZ}}.
\begin{corollary}\label{cor:MertSorZamEF}  {\rm (\cite[Corollary 5]{OSCG})}.  If for each $b\in \bb$ the function $a\mapsto c(a,b)$ is inf-compact, 
then
\begin{equation}\label{eq:classicMSZ}
\min_{\pi^\A\in\P(\A)}\sup_{\pi^\bb\in\P^{fs}(\bb)} \uh(c)(\pi^\A,\pi^\bb)=\sup_{\pi^\bb\in\P^{fs}(\bb)}\min_{\pi^\A\in\P(\A)}\uh(c)(\pi^\A,\pi^\bb).
\end{equation}
\end{corollary}
Corollary~\eqref{cor:DeltaEF} implies additional versions of equality \eqref{eq:classicMSZ} with the set $\P^{fs}(\bb)$ in the left and right hand sides of \eqref{eq:classicMSZ} replaced with arbitrary sets $\Delta^\prime (\bb)$ and $\Delta(\bb)$ specified in statements (a) and (b) of Corollary~\ref{cor:DeltaEF} respectively. In addition, according to \cite[Theorem 17]{OSCG}, equality~\eqref{eq:sup01052019EF} holds when $\bb\subset \Delta^\prime(\bb)\subset \P_{\pi^\A}^S(\bb),$ where  $\P_{\pi^\A}^S(\bb)=\{\pi_*^\bb\in\P(\bb):\, \uh(c)(\pi^\A,\pi_*^\bb)\ \textrm{is well-defined} \},$ $\pi^\A\in\P(\A).$ Therefore, in the left hand side of
\eqref{eq:classicMSZ} the set $\P^{fs}(\bb)$ can be replaced with a set $\Delta^\prime$ such that $ \bb\subset \Delta^\prime(\bb)\subset \P_{\pi^\A}^S(\bb).$

Next we define  solutions for a game with a payoff function unbounded from above and from below.  

\begin{definition}\label{defi:equil}{\rm (Solution for a game).
A pair of   strategies $(\pi^\A,\pi^\bb)\in \P^S(\A)\times\P^S(\bb)$ for Players~I and II is called a \textit{solution  (saddle point, equilibrium) of the game $\{\A,\bb, c\}$} if
\begin{equation}\label{eq:solution}
\uh(c)(\pi^\A,\pi^\bb_*)\le \uh(c)(\pi^\A,\pi^\bb)\le\uh(c)(\pi^\A_*,\pi^\bb)
\end{equation}
for all $\pi_*^\A\in \P(\A)$ and $\pi^\bb_*\in\P(\bb).$ A solution $(\pi^\A,\pi^\bb)$ for the game $\{\A,\bb, c\}$ is  called \textit{pure} if the strategies
$\pi^\A$ and $\pi^\bb$ are pure.}
\end{definition}

If a solution $(\pi^\A,\pi^\bb)\in \P^S(\A)\times\P^S(\bb)$ for the game $\{\A,\bb, c\}$  exists, 
 the  real number $v=\uh(c)^{\flat}(\pi^\bb)=v^\flat =v^\sharp =\uh(c)^{\sharp}(\pi^\A)=\uh(c)(\pi^\A,\pi^\bb)$ is the value of this game. Indeed, since $(\pi^\A,\pi^\bb)\in \P^S(\A)\times\P^S(\bb),$ 
inequalities (\ref{eq:solution}) imply that
$\uh(c)^{\flat}(\pi^\bb)=\uh(c)(\pi^\A,\pi^\bb)=\uh(c)^{\sharp}(\pi^\A).$  Therefore,
\begin{equation}\label{eq2.19EF}
v^\flat=\uh(c)^{\flat}(\pi^\bb)=\uh(c)(\pi^\A,\pi^\bb)=\uh(c)^{\sharp}(\pi^\A)=v^\sharp,
\end{equation}
where the first and the last equalities hold
 since $\uh(c)^\flat(\pi^\bb_*)\le \uh(c)^\sharp(\pi^\A)$ and $\uh(c)^\flat(\pi^\bb)\le \uh(c)^\sharp(\pi^\A_*)$ for all $(\pi^\A_*,\pi^\bb_*)\in \P^S(\A)\times\P^S(\bb).$

Moreover, $(\pi^\A,\pi^\bb)\in \P^S(\A)\times\P^S(\bb)$ is the solution for the game $\{\A,\bb, c\}$ if and only if there exists the value $v,$ the strategy $\pi^\A$ is a solution for Player~I,  and the strategy $\pi^\bb$ is a solution for Player~II. Indeed, the necessary condition is proved in the previous paragraph, and the sufficient condition follows from  $\uh(c)(\pi^\A,\pi^\bb)\le \uh(c)^{\sharp}(\pi^\A)=v=\uh(c)^{\flat}(\pi^\bb)\le
\uh(c)(\pi^\A,\pi^\bb),$ where the equality in the middle hold because $\pi^\A$ and $\pi^\bb$ are solutions for Players I and II respectively, and $v$ is the value.

If  a pure solution $(a,b)\in \A\times\bb$ for the game $\{\A,\bb, c\}$  exists, then the number
\begin{equation}\label{eq:111111}
v=\uh(c)^{\flat}(b)=v^\flat =v^\sharp =\uh(c)^{\sharp}(a)=c(a,b)
\end{equation}
 is the value of this game. Moreover,
\[
v=\sup_{b^*\in\bb}\inf_{a^*\in\A}c(a^*,b^*)=\inf_{a^*\in\A}\sup_{b^*\in\bb}c(a^*,b^*)
\]
because, by the definitions of $\uh(c)^{\flat}$ and $\uh(c)^{\sharp}$,
\[
\sup_{b^*\in\bb}\inf_{a^*\in\A}c(a^*,b^*)=\sup_{b^*\in\bb}\uh(c)^{\flat}(b^*),\quad \inf_{a^*\in\A}\sup_{b^*\in\bb}c(a^*,b^*)=
\inf_{a^*\in\A}\uh(c)^{\sharp}(a^*),\]
and, in view of (\ref{eq:111111}), 
\[
\inf_{a^*\in\A}\uh(c)^{\sharp}(a^*)\le \uh(c)^{\sharp}(a)=v=\uh(c)^{\flat}(b)
\le
\sup_{b^*\in\bb}\uh(c)^{\flat}(b^*)\le \inf_{a^*\in\A}\uh(c)^{\sharp}(a^*),
\]
where the last inequality follows from (\ref{eq:g}) with $\pi^\A=\delta_{\{a\}}$ and $\pi^\bb=\delta_{\{b\}}.$ Therefore, the game $\{\A,\bb, c\}$ has a solution in pure strategies  (that is, the players can play only pure strategies, and this game has a solution)  if and only if there is a pure solution for the game $\{\A,\bb, c\}$.

\begin{remark}\label{rem:sim}
{\rm Let $\{\A,\bb, c\}$ be a two-person zero-sum game introduced in Definition~\ref{defi:game}. Then the triplet $\{\bb,\A,-c^{\A\leftrightarrow\bb}\},$ where $c^{\A\leftrightarrow\bb}(b,a):=c(a,b)$ for each $a\in\A$ and $b\in\bb,$ is also a game satisfying conditions~(i--v) from Definition~\ref{defi:game}.
If this construction is repeated, it leads to the original game.  The game $\{\A,\bb, c\}$ has a value (solution for Player I, Player II, solution)  if and only if  the game $\{\bb,\A,-c^{\A\leftrightarrow\bb}\}$ has a value (solution for Player II, Player I, solution).}
\end{remark}

\section{Main Results and Discussion}\label{Sec:MR}

 Since this paper deals with symmetrically defined games, each assumption or statement for Player I can be reformulated as the corresponding assumption or statement for Player II. In this paper we provide assumptions and statements that involve both Players or only Player I.  Symmetric assumptions and statements for Player II are provided only if they are used in  explanations or in other statements.

To provide sufficient conditions for the validity of {Assumption~\ref{ass:MAIN-NEW-A}},
let us  define
\begin{equation}\label{eq:new1}
\uh(c)^{\oplus,\sharp}(\pi^\A):=(\uh(c)^\oplus)^\sharp(\pi^\A)=\sup_{b\in\bb}\uh(c)^\oplus(\pi^\A,b),\qquad \pi^\A\in\P(\A).
\end{equation}
We also define the symmetric function for Player II,
\begin{equation}\label{eq:new1b}
\uh(c)^{\ominus,\flat}(\pi^\bb):=(\uh(c)^\ominus)^\flat(\pi^\bb)=\inf_{a\in\A}\uh(c)^\ominus(a,\pi^\bb),\quad \pi^\bb\in\P(\bb).
\end{equation}

%
The following theorem describes the sufficient conditions for {Assumption~\ref{ass:MAIN-NEW-A}}.
\begin{theorem}{\rm({Sufficient} conditions for {Assumption~\ref{ass:MAIN-NEW-A}}).}\label{theorem:newA}   
 {Assumption~\ref{ass:MAIN-NEW-A}} holds
 if at least one of the following assumptions is satisfied:
\begin{itemize}
\item[{\rm(L)}] the function $c$ is bounded below on $\A\times\bb;$
\item[{\rm(U)}] the function $c$ is bounded above on $\A\times\bb;$
\item[{\rm(A1)}] there exist $\gamma_\A\in(0,1),$ $L_\A>0,$ and $b_0\in\bb$ such that for each $a\in\A$
\begin{equation}\label{eq:safe2a}
-L_\A+\gamma_\A \uh(c)^{\sharp}(a)\le c^+(a,b_0);
\end{equation}
\item[{\rm(A2)}] there exist $\gamma_\A\in(0,1),$ $L_\A>0,$ and $\pi^\bb_0\in\P^S(\bb)$ such that
\begin{equation}\label{eq:safe2awww}
-\infty<\int_\bb\min\{0,\uh(c)^{\flat}(b)\}\pi^\bb_0(db),\quad\mbox{and}\quad-L_\A+\gamma_\A \uh(c)^{\sharp}(a)\le \uh(c)^\oplus(a,\pi^\bb_0)
\end{equation}
for each $a\in\A;$
\item[{\rm(A3)}] there exist $\gamma_\A\in(0,1)$ and $M_\A>0$ such that  for each $\pi^\A\in\P(\A)$ and $\pi^\bb\in\P(\bb)$
\begin{equation}\label{eq:safe000wwwnew1}
\gamma_\A  {  \uh(c)^\oplus(\pi^\A,\pi^\bb) } \le \uh(c)^\sharp(\pi^\A)+M_\A;
\end{equation}
\item[{\rm(A4)}] there exists a function $\Psi_\A:\R\to\R$ such that $\Psi_\A(s)<+\infty,$ if  $s<+\infty,$ and for each $\pi^\A\in\P(\A)$ and $\pi^\bb\in\P(\bb)$
    \begin{equation}\label{eq:new2}
    \uh(c)^\oplus(\pi^\A,\pi^\bb)\le \Psi_\A(\uh(c)^\sharp(\pi^\A)){;}
    \end{equation}
{\item[{\rm(A5)}]  if $\uh(c)^\sharp(\pi^\A)<+\infty$  for $\pi^\A\in\P(\A),$ then $\uh(c)^{\oplus,\sharp}(\pi^\A)<+\infty.$}
\end{itemize}
{Moreover, Assumption~(A5) is equivalent to statement \eqref{eq:safe1a}.}
\end{theorem}

{The proof of Theorem~\ref{theorem:newA} is provided in Section~\ref{sec:proofs}.}

 The relations between assumptions~(U), (A1)--({A5}) are described  below in implications
\eqref{eq:new3}. According to Theorem~\ref{theorem:newA}, each of these assumption implies {Assumption~\ref{ass:MAIN-NEW-A}}. However,
assumptions~(U), (L), (A1)--({A5}) are simpler than {Assumption~\ref{ass:MAIN-NEW-A}} and useful for applications.

{Assumptions~(L), (U), (A1)--(A5) are formulated in terms of the primitives of the model. For example, Assumption~(A5) means that for { every} probability $\pi^\A$ on the metric space $\A,$ the inequality
\[
\sup_{b\in\bb}\left\{ \int_\A c^+(a,b)\pi^\A(da)-\int_\A c^-(a,b)\pi^\A(da)\right\}<+\infty
\]
implies
\[
\sup_{b\in\bb} \int_\A c^+(a,b)\pi^\A(da)<+\infty.
\]}

{ Simple Example~\ref{exa:1}, which we provide} for illustrative purposes, describes a two-person zero-sum game $\{\A,\bb,c\}$
with noncompact decision sets and unbounded payoffs satisfying assumption~(A1) and, therefore, this game satisfies {Assumption~\ref{ass:MAIN-NEW-A}}.
It also satisfies Assumptions~\ref{AsMain}(a1,a2,b1,b2) and, as follows from Theorem~\ref{eq:safe1b}(B1), it satisfies {the following symmetric statement to Assumption~\ref{ass:MAIN-NEW-A}.
\begin{Assumption}\label{ass:MAIN-NEW-B}
$\uh(c)^\flat(\pi^\bb)< v^\flat$ for all $\pi^\bb\in\P^U(\bb).$
\end{Assumption}
}

\begin{example}\label{exa:1}
{\rm Let $\A=\bb=\mathbb{R},$ $c(a,b)=a^2-b^2,$ $(a,b)\in \mathbb{R}^2.$ Then the game $\{\A,\bb, c \}$ satisfies assumption~(A1) and, therefore, it satisfies {Assumption~\ref{ass:MAIN-NEW-A}}. Indeed, if we consider arbitrary $b_0\in\mathbb{R},$ $\gamma_\A\in(0,1)$ and set $L_\A:=b_0^2,$ then
\[
-L_\A+\gamma_\A \uh(c)^{\sharp}(a) =\gamma_\A a^2-b_0^2\le a^2-b_0^2=c(a,b_0)\le c^+(a,b_0)
\]
for each $a\in\mathbb{R}$ because $\uh(c)^{\sharp}(a)=a^2$ for each $a\in\mathbb{R}.$
}
\end{example}

The following Theorem~\ref{th:exist_value} provides sufficient conditions for the existence of a value and a solution for a one of  the Players for a   two-person zero-sum game with
possibly noncompact decision sets and unbounded  payoffs. This theorem and  Corollary~\ref{cor:exist_value} also describe the properties of the solution sets under these conditions. In general,  an infinite  game may not have  a value; see Yanovskaya~\cite[p. 527]{Yan} and the references to counterexamples by Ville, by Wald, and by Sion and Wolfe cited there. Therefore, some additional conditions for the existence of a value  and solutions are needed.  The  results available in the literature require among other assumptions that either at least one of the decision sets is compact (Mertens et al.~\cite[Propositions I.1.9 and I.2.2]{MSZ},  Alpern and Gal~\cite{AG}, Prokopovich and Yannelis~\cite{PrYa},   Tian~\cite{Tian1992a}) 
or the payoff function is
convex/concave-like (Ansari et al.~\cite{Ansarietal2000}, Aubin and Ekeland~\cite[Theorem 6.2.7]{ObEk}, Nessah and Tian~\cite{NessahandTian2016},
 Perchet and Vigeral~\cite{Per}, 
 Zeng et al.~\cite{Zengetal2006}). Theorems~\ref{th:exist_value} and \ref{teor:mainonestep} of this paper require neither compactness of one of the decision sets nor convexity/concavity-like properties of the payoff function.   Baye et al.~\cite{Baye1993etal}, Nessah and Tian~\cite{NessahandTian2013}, and Tian~\cite{Tian2015} provide sufficient
 conditions for the existence of pure solutions for games with
possibly noncompact decision sets and without convexity-like assumptions. However,   under the assumptions  of Theorems~\ref{th:exist_value} and \ref{teor:mainonestep}, a game may have neither value nor solution when only pure strategies are played; see Proposition~\ref{propsec7}(a). Therefore, the results of this paper do not follow from the references mentioned in this paragraph.
Note that some of these references provide significant advances to the theory of nonzero-sum games and consider relaxed versions of  semi-continuity assumptions on payoff functions such as transfer lower and upper semi-continuity. However, the results are currently available either for compact decision sets or for convex/concave-like payoff functions or games with solutions in pure strategies.

%
%
%

\begin{theorem}{\rm (Existence of a value and solution for Player~I).} \label{th:exist_value}
Let a  two-person zero-sum game $\{\A,\bb, c \}$   satisfy Assumptions~\ref{AsMain}(a1,a2) and
{\begin{equation}\label{ass:MAIN-NEW-A-WEAK}
\uh(c)^\sharp(\pi^\A)\ge v^\sharp\mbox{ for all }\pi^\A\in\P^U(\A).
\end{equation}}
Then the game $\{\A,\bb, c \}$ has the value $v,$
that is, equality
\eqref{eq:valval} and, therefore, equalities (\ref{eq:valvalclass}) hold.
Moreover, {if Assumption~\ref{ass:MAIN-NEW-A} holds, then}  the set $\P^\sharp_{v}(\A),$  which is a subset of $\P^S(\A),$ is the set of solutions for Player~I, and  $\P^\sharp_{v}(\A)$
is a nonempty convex compact subset of $\P(\A).$
\end{theorem}

{The proof of Theorem~\ref{th:exist_value} is provided in Section~\ref{sec:proofs}.}


\begin{corollary}\label{cor:exist_value}{\rm (Compactness of the set $\P^\flat_{v}(\bb)$).}
Let  the assumptions of Theorem~\ref{th:exist_value} and Assumptions~\ref{AsMain}(b1, b2) hold. Then, in addition to the conclusions of  Theorem~\ref{th:exist_value}, the set $\P^\flat_{v}(\bb)$ is a nonempty convex compact subset of $\P(\bb).$
\end{corollary}
\begin{proof}
 Since all the conditions of Theorem~\ref{th:exist_value} are also assumed in the corollary, the conclusions of Theorem~\ref{th:exist_value} hold. In particular, the game has the value.
The additional conclusions follow from Theorem~\ref{th:exvalPaper1}(a) applied to the game $\{\bb,\A, -c^{\A\leftrightarrow\bb} \}$ introduced in Remark~\ref{rem:sim}.
\end{proof}
\begin{corollary}{\rm (Existence of a value and solution for Player~II).}\label{corth:exist_value}
Let a  two-person zero-sum game $\{\A,\bb, c \}$  
satisfy Assumptions~\ref{AsMain}(b1,b2) and {\begin{equation}\label{ass:MAIN-NEW-B-WEAK}
\uh(c)^\flat(\pi^\bb)\le v^\flat\mbox{ for all }\pi^\bb\in\P^U(\bb).
\end{equation}}
Then the game $\{\A,\bb, c \}$ has the value $v.$
Moreover, {if Assumption~\ref{ass:MAIN-NEW-B} holds, then} the set  $\P^\flat_{v}(\bb),$ which is a subset of $\P^S(\bb),$ is the set of solutions for Player~II, and  $\P^\flat_{v}(\bb)$
is a nonempty convex compact subset of $\P(\bb).$
\end{corollary}
\begin{proof}
The corollary follows from  Theorems~\ref{th:exist_value} applied to the game $\{\bb,\A, -c^{\A\leftrightarrow\bb} \}.$
\end{proof}

 The following statement
\begin{equation}\label{eq:safe1b}
\uh(c)^\flat(\pi^\bb)=-\infty\mbox{ for all } \pi^\bb\in\P^U(\bb),
\end{equation}
is similar and symmetric to statement~\eqref{eq:safe1a}.
Note that statement~(\ref{eq:safe1b}) holds if and only if $\P^\flat_{\alpha}(\bb)\subset \P^S(\bb)$ for each $\alpha \in\mathbb{R}.$
 The following theorem is similar and symmetric to Theorem~\ref{theorem:newA}.
\begin{theorem}{\rm ({Sufficient} conditions for Assumption~{\ref{ass:MAIN-NEW-B}).}} \label{theorem:newB}
 {Assumption~\ref{ass:MAIN-NEW-B}} holds
 if at least one of the following assumptions is satisfied:
\begin{itemize}
\item[{\rm(C)}] either assumption~(L) or assumption (U) holds;
\item[{\rm(B1)}] there exist $\gamma_\bb\in(0,1),$ $L_\bb>0,$ and $a_0\in\A$ such that for each $b\in\bb$
\begin{equation}\label{eq:safe2b}
c^-(a_0,b)\le \gamma_\bb \uh(c)^{\flat}(b)+L_\bb;
\end{equation}
\item[{\rm(B2)}] there exist $\gamma_\bb\in(0,1),$ $L_\bb>0,$ and $\pi^\A_0\in\P(\A)$ such that
\begin{equation}\label{eq:safe2bwww}
\int_\A\max\{0,\uh(c)^{\sharp}(a)\}\pi^\A_0(da)<+\infty,\quad\mbox{and}\quad\uh(c)^\ominus(\pi^\A_0,b)\le \gamma_\bb \uh(c)^{\flat}(b)+L_\bb
\end{equation}
for each $b\in\bb;$
\item[{\rm(B3)}] there exist $\gamma_\bb\in(0,1)$ and $M_\bb>0$ such that  for each $\pi^\A\in\P(\A)$ and $\pi^\bb\in\P(\bb)$
\begin{equation}\label{eq:safe000bbwwwnew1}
\uh(c)^\flat(\pi^\bb)\le \gamma_\bb {  \uh(c)^\ominus(\pi^\A,\pi^\bb) }
+M_\bb;
\end{equation}
\item[{\rm(B4)}] there exists a function $\Psi_\bb:\R\to\R$ such that $\Psi_\bb(s)>-\infty,$ if  $s>-\infty,$ and for each $\pi^\A\in\P(\A)$ and $\pi^\bb\in\P(\bb)$
    \begin{equation}\label{eq:new2b}
    \Psi_\bb(\uh(c)^\flat(\pi^\bb))\le\uh(c)^\ominus(\pi^\A,\pi^\bb){;}
    \end{equation}
{\item[{\rm(B5)}]  if $\uh(c)^\flat(\pi^\bb)>-\infty$ for $\pi^\bb\in\P(\bb),$ then $\uh(c)^{\ominus,\flat}(\pi^\bb)>-\infty.$}
\end{itemize}
{Moreover, Assumption~(B5) is equivalent to statement \eqref{eq:safe1b}.}
\end{theorem}

{The proof of Theorem~\ref{theorem:newB} is provided in Section~\ref{sec:proofs}.}

\begin{remark}\label{rem:new}
{\rm   Assumptions~(C), (B1)--{(B5)} are useful for applications. A two-person zero-sum game $\{\A,\bb, c \}$ 
satisfies  {A}ssumption~(L) ({Assumptions~\ref{ass:MAIN-NEW-B}}, (U){,} (C), (B1), (B2), (B3), (B4){, (B5)} respectively) if and only if
the game $\{\bb,\A,-c^{\A\leftrightarrow\bb}\}$ introduced in Remark~\ref{rem:sim} satisfies {A}ssumptions~(U) ({Assumptions~\ref{ass:MAIN-NEW-A}}, (L), (C), (A1), (A2), (A3), (A4){, (A5)} respectively).}
\end{remark}

The following theorem describes sufficient conditions for the existence of a solution for a game  and the structure of the solution set.

\begin{theorem}\label{teor:mainonestep}{\rm (Existence of a solution for a game).}
Let a two-person zero-sum game $\{\A,\bb, c \}$ 
satisfy Assumptions~\ref{AsMain}(a1,a2,b1,b2) and Assumptions~{\ref{ass:MAIN-NEW-A}, \ref{ass:MAIN-NEW-B}}.  Then:
\begin{itemize}
\item[\rm(i)] the game $\{\A,\bb, c \}$ has a value $v \in\mathbb{R}$ and a solution $(\pi^\A,\pi^\bb)\in \P^\sharp_{v }(\A)\times \P_{v }^\flat(\bb);$
\item[\rm(ii)]  the sets  $\P^\sharp_{v }(\A)$ and $\P_{v }^\flat(\bb)$
 are nonempty convex compact subsets of $\P(\A)$ and $\P(\bb)$  respectively; moreover, $\P^\sharp_{v }(\A)\subset\P^S(\A)$ and $\P_{v }^\flat(\bb)\subset\P^S(\bb);$
\item[\rm(iii)] a pair of
  strategies $(\pi^\A,\pi^\bb)\in \P^S(\A)\times \P^S(\bb)$
  is a solution for the game $\{\A,\bb, c \}$ if and only if
  $\pi^\A\in \P^\sharp_{v }(\A)$ and  $\pi^\bb\in \P_{v }^\flat(\bb).$
  \end{itemize}
\end{theorem}

{The proof of Theorem~\ref{teor:mainonestep} is provided in Section~\ref{sec:proofs}.}

\begin{remark}\label{rem:cond(iii,iv)}
{\rm As explained in Remark~\ref{rem:cond(i,ii)}, inf-compactness in $a$ of the function $c(a,b)$ stated in Assumption~\ref{AsMain}(a2)  is close to the necessary condition for the existence of a solution for Player~I.  If the function $c$ is bounded above on $\A\times\bb,$ then the set $\A$ is compact.  The similar observation takes place for sup-compactness in $b$ of the function $c(a,b)$ stated in Assumption~\ref{AsMain}(b2), the existence of a solution for Player~II, and compactness of the set $\bb.$  We observe that the expected payoff $\uh(c)(\pi^\A,\pi^\bb)$ is well-defined for all pairs of strategies $(\pi^\A,\pi^\bb)\in \P(\A)\times \P(\bb)$ if and only if the function $c(a,b)$ is bounded either from above or from below on $\A\times\bb.$  In these two cases, the sets $\A$ and $\bb$ are compact respectively.  If the function $c$ is bounded on $\A\times\bb,$ then under  Assumptions~\ref{AsMain}(a2,b2)  both decision sets are compact.  So, for a problem with two noncompact decision sets $\A$ and $\bb,$ under natural inf/sup-compactness Assumptions~\ref{AsMain}(a2,b2) 
the payoff function is  unbounded from above and from below on $\A\times\bb,$ and there exist pairs of strategies $(\pi^\A,\pi^\bb)\in \P(\A)\times \P(\bb)$ with undefined values of $\uh(c)(\pi^\A,\pi^\bb).$}
\end{remark}

If the function $c$ is bounded above  or below on $\A\times\bb,$ then Assumptions~{\ref{ass:MAIN-NEW-A} and \ref{ass:MAIN-NEW-B}}  hold; see Theorems~\ref{theorem:newA} and \ref{theorem:newB}.
As follows from these observations and Remark~\ref{rem:cond(i,ii)}, Theorems~\ref{th:exist_value}, \ref{teor:mainonestep} and Corollary~\ref{cor:exist_value} imply several known results for games with  two compact decision sets and with at least one compact decision set.
In particular, Theorems~\ref{th:exist_value}  and Corollary~\ref{corth:exist_value} imply Glicksberg's theorem: for a game with two compact decision sets the value exists if the payoff function is upper (or lower) semi-continuous.  If $\A$ and $\bb$ are subsets of Polish spaces, then  Theorem~\ref{teor:mainonestep} implies \cite[Theorem I.2.4]{MSZ} stating that the game has a solution if $\A$ and $\bb$ are compact sets, $c$ is a real-valued function bounded from below or from above, and Assumptions~\ref{AsMain}(a1,b1) hold (that is, for each $b\in \bb$ the function $a\mapsto c(a,b)$ is lower semi-continuous, and for each $a\in \A$ the function $b\mapsto c(a,b)$ is upper semi-continuous). Indeed, if $c$ ia a real-valued function bounded from below or bounded from above, then, in view of   Theorem~\ref{theorem:newA}(L, U) and Theorem~\ref{theorem:newB}(C), Assumptions~{\ref{ass:MAIN-NEW-A} and \ref{ass:MAIN-NEW-B}} hold, and, if the sets $\A$ and $\bb$ are compact, then Assumption~\ref{AsMain}(a1) implies Assumption~\ref{AsMain}(a2), and  Assumption~\ref{AsMain}(b1) implies Assumption~\ref{AsMain}(b2). We also remark that \cite[Theorem I.2.4]{MSZ} is a more general statement that the version of von Neumann's theorem for mixed strategies, which states the existence of a solution, if the  decision sets are compact and the payoff function is continuous; Owen~\cite[Theorem IV.6.1]{Owen} or Petrosyan and Zenkevich~\cite[Theorem 2.4.4]{PZ}.

 The following corollary from Theorem~\ref{teor:mainonestep} generalizes \cite[Theorem 20]{OSCG}.
\begin{corollary} {\rm (Existence of a solution for  a game)}\label{cor:new}
Let a  two-person zero-sum game $\{\A,\bb, c \}$  
satisfy Assumptions~\ref{AsMain}(a1,a2,b1,b2).  If, in addition, the payoff function $c$ is bounded  either below or above on $\A\times\bb,$ then the conclusions of Theorem~\ref{teor:mainonestep} hold.
\end{corollary}

{The proof of Corollary~\ref{cor:new} is provided in Section~\ref{sec:proofs}.}

 As explained in Remark~\ref{rem:cond(i,ii)}, if the function $c$ is bounded above (below) in Corollary~\ref{cor:new}, then the set $\A$ ($\bb$) is compact.   The following example  demonstrates that Assumption~\ref{AsMain}(b2) is essential in Corollary~\ref{cor:new}  when the function $c$ is bounded below on $\A\times\bb$. Of course, this is also true for  Assumption~\ref{AsMain}(a2) when the function $c$ is bounded above on $\A\times\bb$.
\begin{example}\label{exa:2}
{\rm   This example describes
a two-person zero-sum game $\{\A,\bb, c \}$  
 with the payoff function $c$  bounded from below on $\A\times\bb$ and satisfying Assumptions~\ref{AsMain}(a1,a2,b1).   However,  the function $b\mapsto c(a,b)$ is not sup-compact on $\bb$ for each $a\in \A,$ and $\P_{v }^\flat(\bb)=\emptyset.$ Therefore, this game has no solution.

Let $\A=\bb:=\mathbb{R},$  and 
\[
c(a,b):=1+a^2-\frac{\exp(b)}{1+\exp(b)},\quad a,b\in\mathbb{R}.
\]
Note that the function $c$ takes positive values and it is continuous on $\mathbb{R}^2.$ Moreover, the function $a\mapsto c(a,b)$ is obviously inf-compact on $\mathbb{R}$ for each $b\in\mathbb{R}.$  However, for
each $a\in\mathbb{R}$ the function $b\mapsto c(a,b)$ is not sup-compact on $\mathbb{R}$ because for every $a\in\mathbb{R}$ the set $\{b\in \mathbb{R}\,:\, c(a,b)\ge 0\}=\mathbb{R}$ is not compact.

The set $\P_{v }^\flat(\bb)$ is empty. Indeed, direct calculations imply that $v^\sharp=v^\flat=0$ and $\uh(c)^\flat(\pi^\bb)=1-\int_{\mathbb{R}}\frac{\exp(b)}{1+\exp(b)}\pi^\bb(db)$ for each
$\pi^\bb\in\P(\mathbb{R}).$ Therefore,
\[
\P_{v }^\flat(\bb)=\{\pi^\bb\in\P(\mathbb{R})\,:\, \int_{\mathbb{R}}\frac{\exp(b)}{1+\exp(b)}\pi^\bb(db)=1\}=\emptyset,
\]
where the last equality holds because $\pi^\bb(\mathbb{R})=1$ and $\frac{\exp(b)}{1+\exp(b)}<1$ for each $b\in\mathbb{R}.$
The game $\{\A,\bb, c \}$ has no solution
since the set $\P_{v }^\flat(\bb)$ is empty.
}
\end{example}

\section{Number Guessing Game}\label{sec:ex}

In this section we consider the following game to illustrate the results of this paper.
Two players select nonnegative numbers $a$ and $b,$ and Player I pays the amount of
$c(a,b)=\varphi(a-b)$ to Player II.  For example, if the player, who selects the larger number wins, that is,
$\varphi(a-b)={\rm I}(a>b),$  this game does not have a value; see e.g., \cite{Yan}. 
We apply the results of our paper to such games.  In particular, for a polynomial function $\varphi,$   Proposition~\ref {propsec7} completely characterizes  all the situations when the games have values and solutions.

 Since  both decision sets $\A=\bb:=\mathbb{R}_+=[0,+\infty)$ are not compact, the only previously available result on the existence of the solution is Aubin and Ekeland~\cite[Theorem 6.2.7]{ObEk}, which assumes Assumptions~\ref{AsMain}(a1,a2), the concavity of $c(a,b)$ in $a$ and convexity  of $c(a,b)$ in $b.$  Under these conditions, there exists a pure solution  for the game.    Another simple sufficient condition, under which a two-person zero-sum game $\{\A,\bb, c \}$ 
 has a pure solution, is  $\A=\bb:=\mathbb{R}_+$ and the function $c(a,b)$ is nondecreasing in $a$ and nonincreasing  in $b.$  In this case it is optimal for each player to select the decision $0.$ These arguments are applicable to Example~\ref{exa:1}.

 The examples provided in this section may satisfy neither of the two described sufficient conditions.  In addition, according to Proposition~\ref{propsec7},  solutions in pure strategies may not exist for the provided examples.

 \begin{example}\label{exasec7A}
 {\rm Let
$\A=\bb:=\mathbb{R}_+$ and $c(a,b):=\varphi(a-b)$ for each $a,b\in \mathbb{R}_+,$ where $\varphi:\mathbb{R}\to\mathbb{R}$ is a continuous function. To be consistent with assumptions (iv, v) in Definition~\ref{defi:game}, assume that
\begin{equation}\label{eq:assgame}
-\infty< \ilim_{s\to+\infty} \varphi(s) \quad \mbox{and}\quad \slim_{s\to-\infty} \varphi(s)<+\infty.
\end{equation}

The triple  $\{\A,\bb, c \}$ is a two-person zero-sum game introduced in Definition~\ref{defi:game} because the function $a\mapsto c(a,b)$ is bounded  below on $\A$ for each $b\in\bb,$ if the first inequality in (\ref{eq:assgame}) holds, and the function $b\mapsto c(a,b)$ is bounded  above on $\bb$
for each $a\in\A,$ if the second inequality in (\ref{eq:assgame}) holds.
}
\end{example}

\begin{proposition}\label{propsec7A}  Consider the two-person zero-sum game defined in Example~\ref{exasec7A}. Then:
  \begin{itemize}
  \item[(a)] if $\varphi(s)\to+\infty$ as $s\to+\infty,$ then  Assumptions~\ref{AsMain}(a1,a2) and therefore the conclusions of Theorem~\ref{th:exvalPaper1} hold;
  \item[(b)] if $\varphi(s)\to+\infty$ as $s\to+\infty$ and $\varphi(s)=\varphi_1(s)+\varphi_2(s)$ for each $s\in\mathbb{R},$ where $\varphi_1:\mathbb{R}\to\mathbb{R}$ is increasing and $\varphi_2:\mathbb{R}\to\mathbb{R}$ is bounded, then  the assumptions and therefore the conclusions of Theorem~\ref{th:exist_value} hold;
  \item[(c)] if $\varphi(s)\to+\infty$ as $s\to+\infty,$ $\varphi(s)\to-\infty$ as $s\to-\infty,$ and $\varphi(s)=\varphi_1(s)+\varphi_2(s)$ for each $s\in\mathbb{R},$
where $\varphi_1:\mathbb{R}\to\mathbb{R}$ is increasing and $\varphi_2:\mathbb{R}\to\mathbb{R}$ is bounded, then the assumptions and therefore the conclusions of Theorem~\ref{teor:mainonestep} hold.
  \end{itemize}
\end{proposition}

{The proof of Proposition~\ref{propsec7A} is provided in Section~\ref{sec:proofs}.}

\begin{example}\label{exasec7}
{\rm Consider Example~\ref{exasec7A} with the function $\varphi$ being
 a polynomial of a degree $M=1,2,\ldots,$ that is, $\varphi(s)=\sum_{n=0}^M\alpha_ns^n,$ $s\in\mathbb{R},$ where $\alpha_n\in\mathbb{R},$ $n=0,\ldots,M,$ and $a_M\ne 0.$}
\end{example}

\begin{proposition}\label{propsec7}  Consider the two-person zero-sum game defined in Example~\ref{exasec7}.

  (a) If the  integer $M$ is odd and $\alpha_M>0,$ then  this game satisfies the assumptions of Theorem~\ref{teor:mainonestep}, and therefore the conclusions of Theorem~\ref{teor:mainonestep}  hold  for this game.  Furthermore, if $M\ge 3$ and $\alpha_1<0,$ then there is no pure solution for this game.

  (b) If the integer $M$ is even or $\alpha_M<0,$ then  this game does not satisfy either assumption (iv) or assumption (v) from Definition~\ref{defi:game}, and there is no finite value because either $|v^\flat|=+\infty$ or $|v^\sharp|=+\infty.$
\end{proposition}

{The proof of Proposition~\ref{propsec7} is provided in Section~\ref{sec:proofs}.}

We note that in the two-person zero-sum game from Proposition~\ref{propsec7}(a) action $0$ for each Player strongly dominates any other action large enough and after elimination of these actions we have a compact game. On the other hand, let us consider the two-person zero-sum game defined in Example~\ref{exasec7A}. If the assumptions of Proposition~\ref{propsec7A}(c) hold, and for each $s^*\in\mathbb{R}\setminus\{0\}$
\begin{equation}\label{eq:nnew1}
\ilim_{s\to\pm\infty} (\varphi(s)-\varphi(s+s^*))<0<\slim_{s\to\pm\infty} (\varphi(s)-\varphi(s+s^*)),
\end{equation}
then an any action of each Player does not dominates any other his/her action because, according to \eqref{eq:nnew1}, for each $a_*,a^*\in \mathbb{R}_+,$ $a_*\ne a^*,$ there exist $b_*,b^*\in\mathbb{R}_+$ such that
\[
c(a_*,b_*)-c(a^*,b_*)=\varphi(a_*-b_*)-\varphi(a^*-b_*)<0<\varphi(a_*-b^*)-\varphi(a^*-b^*)=c(a_*,b^*)-c(a^*,b^*),
\]
and, symmetrically,
for each $b_*,b^*\in \mathbb{R}_+,$ $b_*\ne b^*,$ there exist $a_*,a^*\in\mathbb{R}_+$ such that
\[
c(a^*,b_*)-c(a_*,b_*)=\varphi(a^*-b_*)-\varphi(a_*-b_*)<0<\varphi(a_*-b_*)-\varphi(a_*-b^*)=c(a_*,b_*)-c(a_*,b^*).
\]
The example of such function is $\varphi=\varphi_1+\varphi_2$ with $\varphi_1(s)=${\rm sgn}$(s)\ln(|s|+1)$ and ${\varphi_2(s)={\sin(s)+\sin(\sqrt2 s)}},$ $s\in\mathbb{R}.$ Indeed, the assumptions of Proposition~\ref{propsec7A}(c) are trivial,
inequalities \eqref{eq:nnew1} hold because $\varphi_1(s)-\varphi_1(s+s^*)\to 0$ as $s\to\pm\infty,$ and $\varphi_2$ satisfies \eqref{eq:nnew1}
for each $s^*\in\mathbb{R}\setminus\{0\}$ since the functions $s\mapsto \sin(s)$ and $s\mapsto (\sqrt2 s)$ have commensurable prime periods.

We notice that \cite[Theorem 6.2.7]{ObEk} {cannot be applied}  in most cases   to the examples considered in this section because it assumes   concavity of $c(a,b)$ in $a$ and convexity  of $c(a,b)$ in $b.$ For example, assume that the function $\varphi$ is twice differentiable, as this holds in Example~\ref{exasec7}.  Then $\frac{\partial^2 \varphi(a-b)}{\partial a^2}=\frac{\partial^2 \varphi(a-b)}{\partial b^2}.$  Therefore, the convexity/concavity assumption implies that these derivatives are equal to 0, and $\varphi(a-b)=M(a-b)+C.$ For $M>0$ this game is covered by Proposition~\ref{propsec7}(a), and $a=b=0$ is the solution.     For $M<0$ this game is covered by Proposition~\ref{propsec7}(ii), and there is no  solution.

We also remark that the last claim of Proposition~\ref{propsec7}(a), which states nonexistence of a pure solution, in some sense complements the result by Dreshen, Karlin, and Shapley (see Parrilo~\cite[Theorem 2.2]{PP}) that states that, for a  game with $\A=\bb=[0,1]$ and with a polynomial payoff function $c,$ there exists a solution $(\pi^\A,\pi^\bb)$ with $\pi^\A$ and $\pi^\bb$ having finite supports. In the case of a polynomial function $\varphi$ defined in Example~\ref{exasec7}, each of these finite supports  consists of no more than $(M+1)$ points.

\section{Proofs}\label{sec:proofs}

This section consists of four subsections. Subsection~\ref{newsubsec:1} provides  the proofs of Theorems~\ref{theorem:newA}, \ref{th:exist_value}, and \ref{theorem:newB}, 
Subsection~\ref{subsec:2} provides the proofs of Theorem~\ref{teor:mainonestep}  and Corollary~\ref{cor:new}, and Subsection~\ref{subsec:3} provides the proofs of Propositions~\ref{propsec7A} and \ref{propsec7}.

\subsection{Proofs of Theorems~\ref{theorem:newA}, \ref{th:exist_value}, and \ref{theorem:newB}}\label{newsubsec:1}

The proof of  Theorem~\ref{theorem:newA} 
is based on  Lemmas~{\ref{lem:safe1bbb27}--\ref{lem:safe1bbbbbbaaaaanew}.}

\begin{lemma}\label{lem:safe1bbb27}
Consider  a two-person zero-sum game $\{\A,\bb, c \}.$ 
If $\pi^\A\in\P^S(\A),$ then either $\uh(c)^{\oplus,\sharp}(\pi^\A)<+\infty$ or
$\inf_{b\in\bb} \uh(c)^{\ominus}(\pi^\A,b)>-\infty.$
\end{lemma}
\begin{proof}
On the contrary, let   $\uh(c)^{\oplus,\sharp}(\pi^\A)=+\infty$ and    $\inf_{b\in\bb} \uh(c)^{\ominus}(\pi^\A,b)=-\infty$ for a strategy $\pi^\A\in\P^S(\A).$  Then for each $n=1,2,\ldots$ there exist two points $b^{(1)}_n, b^{(2)}_n\in \bb$ such that $\uh(c)^{\oplus}(\pi^\A,b^{(1)}_n)\ge 2^{n}$ and $\uh(c)^{\ominus}(\pi^\A,b^{(2)}_n)\le -2^{n}.$  Let us consider the probability measures $\pi_i^\bb(B)=\sum_{n=1}^\infty 2^{-n}I\{b^{(i)}_n\in B\}$ for $B\in{\mathcal B}(\bb),$   $i=1,2.$   We define $\pi^\bb=\frac{1}{2}(\pi_1^\bb+\pi_2^\bb).$  Then, by Fubini's theorem, $\uh(c)^{\oplus}(\pi^\A,\pi_1^\bb)=\sum_{n=1}^\infty \pi(b^{(1)}_n) \uh(c)^{\oplus}(\pi^\A,b^{(1)}_n)\ge\sum_{n=1}^\infty 2^{-n}2^n=+\infty.$  Therefore,
$\uh(c)^{\oplus}(\pi^\A,\pi_1^\bb)=+\infty$, and $\uh(c)^{\oplus}(\pi^\A,\pi^\bb)= \frac{1}{2}\uh(c)^{\oplus}(\pi^\A,\pi_1^\bb)+ \frac{1}{2}\uh(c)^{\oplus}(\pi^\A,\pi_2^\bb)= +\infty.$ Similarly, $\uh(c)^{\ominus}(\pi^\A,\pi_2^\bb)=-\infty$ which implies $\uh(c)^{\ominus}(\pi^\A,\pi^\bb)=-\infty.$  Thus, $\pi^\A\notin\P^S(\A).$
\end{proof}
\begin{lemma}\label{lem:safe1bbbbbbaaaaa}
{Assumption~(A5) of Theorem~\ref{theorem:newA}} is equivalent  to statement \eqref{eq:safe1a}.
\end{lemma}
\begin{proof}Statement \eqref{eq:safe1a} is equivalent  to its  contrapositive statement
\begin{equation}\label{eq:equiv lem}
\mbox{if }\uh(c)^\sharp(\pi^\A)<+\infty \mbox{  for }  \pi^\A\in\P(\A),\mbox{ then } \pi^\A\in\P^S(\A).
\end{equation}

Let {Assumption~(A5)} hold.
Suppose $\pi^\A\in\P(\A)$ satisfies the inequality $\uh(c)^\sharp(\pi^\A)<+\infty.$ Observe that $\pi^\A\in\P^S(\A).$
Indeed, since $\uh(c)^\sharp(\pi^\A)<+\infty,$ then {Assumption~(A5)} implies 
$
 {  \uh(c)^\oplus(\pi^\A,\pi^\bb) } <+\infty
$
for each $\pi^\bb\in\P(\bb),$ that is,  $\pi^\A\in\P^S(\A).$  Thus, statement \eqref{eq:equiv lem} holds.

Now let  statement \eqref{eq:equiv lem} hold. We consider an arbitrary $\pi^\A\in\P(\A)$ satisfying the inequality  $\uh(c)^\sharp(\pi^\A)<+\infty.$ Then, in view of \eqref{eq:equiv lem}, $\pi^\A\in\P^S(\A).$  Lemma~\ref{lem:safe1bbb27} implies that either $\uh(c)^{\oplus,\sharp}(\pi^\A)<+\infty$ or
$\inf_{b\in\bb} \uh(c)^{\ominus}(\pi^\A,b)>-\infty.$  To complete the proof of the validity of {Assumption~(A5)}, it is sufficient to show that the latter inequality implies the former one.  Indeed, let $\inf_{b\in\bb}\uh(c)^{\ominus}(\pi^\A,b)>-\infty.$ Therefore, $\uh(c)^{\oplus,\sharp}(\pi^\A)=\sup_{b\in\bb}\{\uh(c)(\pi^\A,b)-\uh(c)^{\ominus}(\pi^\A,b)\}\le \uh(c)^\sharp (\pi^\A)-\inf_{b\in\bb} \uh(c)^{\ominus}(\pi^\A,b)<+\infty.$ 
\end{proof}

\begin{lemma}\label{lem:safe1bbbbbbaaaaanew}  For a two-person zero-sum game
 $\{\A,\bb, c \},$     the following implications hold for the assumptions introduced in Theorem~\ref{theorem:newA}:
\begin{equation}\label{eq:new3}
{\rm(U)}\Rightarrow{\rm(A1)}\Rightarrow{\rm(A2)}\Rightarrow{\rm(A3)}\Rightarrow{\rm(A4)}{\Rightarrow{\rm(A5)}}.
\end{equation}
\end{lemma}

\begin{proof}
``${\rm(U)}\Rightarrow{\rm(A1)}$'': If the function $c$ is bounded above on $\A\times\bb,$ then inequality \eqref{eq:safe2a} takes place for  $L_\A=\max\{0,\sup\{c(a,b):a\in\A,b\in\bb\}\}\in[0,+\infty),$   for all $ \gamma_\A\in(0,1),$  and   for all $b_0\in\bb.$ 

``${\rm(A1)}\Rightarrow{\rm(A2)}$'': assumption~(A1) implies assumption~(A2)  with $\pi^\bb_0=\delta_{\{b_0\}}$  since \eqref{eq:safe2a} becomes the second inequality in \eqref{eq:safe2awww}, and the first one becomes   $\uh(c)^{\flat}(b_0)>-\infty,$ which is true in view of
 (\ref{eq:safe00}). 

``${\rm(A2)}\Rightarrow{\rm(A3)}$'':
Let us fix arbitrary $\pi^\A\in\P(\A),\,\pi^\bb\in\P(\bb)$ and prove  that
\begin{equation}\label{eq:safe000www}
\gamma_\A  {  \uh(c)^\oplus(\pi^\A,\pi^\bb) } 
\le \uh(c)^\sharp(\pi^\A)+L_\A-\int_{\bb}\min\{0,\uh(c)^{\flat}(b)\}\pi^\bb_0(db).
\end{equation} Note that
\begin{equation}\label{eq:safe001}
\sup_{b\in\bb}c^+(a,b)=\max\{\uh(c)^{\sharp}(a),0\}   { ,\qquad\qquad a\in\A.}
\end{equation}
Indeed, if $\uh(c)^{\sharp}(a)\le0,$ then $c(a,b)\le 0$ for each $b\in\bb,$  and  both sides of \eqref{eq:safe001} equal $0.$ If $\uh(c)^{\sharp}(a)>0,$ then the set $B^+(a):=\{b\in\bb\,:\, c(a,b)>0\}=\{b\in\bb\,:\, c^+(a,b)>0\}$ is nonempty, and for each $a\in\A$
\[
\sup_{b\in \bb}c^+(a,b)=\sup_{b\in B^+(a)}c^+(a,b)=\sup_{b\in B^+(a)}c(a,b)=\sup_{b\in \bb}c(a,b)=\uh(c)^{\sharp}(a) =\max\{\uh(c)^{\sharp}(a),0\};
\]
{these equalities} follow from the basic properties of suprema and the definition of $B^+(a).$

Equality (\ref{eq:safe001}) implies that for all $\pi^\A\in\P(\A)$ and for all $\pi^\bb\in\P(\bb)$
\begin{equation}\label{eq:safe002}
\begin{aligned}
& {  \uh(c)^\oplus(\pi^\A,\pi^\bb) } = \int_\A\int_\bb c^+(a,b)\pi^\bb(db)\pi^\A(da)\\
& \le
{ \int_\A\int_\bb\max\{\uh(c)^{\sharp}(a),0\}\pi^\bb(db)\pi^\A(da)= }
\int_\A\max\{\uh(c)^{\sharp}(a),0\}\pi^\A(da).
\end{aligned}
\end{equation}

The  second inequality in (\ref{eq:safe2awww}) implies 
\begin{equation}\label{epreintEF}
\gamma_\A \max\{\uh(c)^{\sharp}(a),0\}\le \uh(c)^\oplus(a,\pi^\bb_0)+L_\A {,\qquad\qquad a\in\A.}
\end{equation}

The integration of  both sides of \eqref{epreintEF} in $\pi^\A \in \P(\A)$ leads to
\begin{equation}\label{eq:safe003aa}
\gamma_\A \int_\A\max\{\uh(c)^{\sharp}(a),0\}\pi^\A(da)\le \uh(c)^\oplus(\pi^\A,\pi^\bb_0)+L_\A  =\uh(c)(\pi^\A,\pi_0^\bb)-\uh(c)^\ominus(\pi^\A,\pi^\bb_0)+L_\A.
\end{equation}
 Observe that
\begin{equation}\label{eq:safe003aaa}
-\infty<\int_\bb\min\{0,\uh(c)^{\flat}(b)\}\pi^\bb_0(db) {= \int_\A\int_\bb\min\{0,\uh(c)^{\flat}(b)\}\pi^\bb_0(db)\pi^\A(da) }
\le\uh(c)^\ominus(\pi^\A,\pi^\bb_0),
\end{equation}
where the first inequality in (\ref{eq:safe003aaa}) {is} the first inequality in (\ref{eq:safe2awww}), {the equality follows from integrating the constant in $\pi^\A,$ } and the last inequality
{follows from the second inequality in} (\ref{eq:safe00}).

Inequalities (\ref{eq:safe003aa}) and (\ref{eq:safe003aaa})  imply that
\begin{equation}\label{eq:safe003}
\gamma_\A \int_\A\max\{\uh(c)^{\sharp}(a),0\}\pi^\A(da)\le\uh(c)(\pi^\A,\pi^\bb_0)+L_\A-\int_\bb\min\{0,\uh(c)^{\flat}(b)\}\pi^\bb_0(db).
\end{equation}
{Observe that $-\infty<\uh(c)(\pi^\A,\pi^\bb_0)\le
\uh(c)^\sharp(\pi^\A)$   because of (\ref{eq:safe003aaa}) and the definition of $\uh(c)^\sharp(\pi^\A).$  Therefore,   inequalities (\ref{eq:safe002}) and (\ref{eq:safe003}) imply (\ref{eq:safe000www}).
}

 Let us  choose $\gamma_\A\in(0,1)$ and $M_\A:=\max\{L_\A-\int_\bb\min\{0,\uh(c)^{\flat}(b)\}\pi^\bb_0(db),1\}>0.$ Note that the first inequality in (\ref{eq:safe2awww}) implies $M_\A<+\infty.$ Thus (\ref{eq:safe000wwwnew1}) follows from (\ref{eq:safe000www}).

``${\rm(A3)}\Rightarrow{\rm(A4)}$'':  Inequality (\ref{eq:new2}) follows from \eqref{eq:safe000wwwnew1} if we set $\Psi_\A(s):=\frac{1}{\gamma_\A}(s+{M_\A})$ for each $s\in\R.$

{``${\rm (A4)}\Rightarrow{\rm (A5)}$'': Inequality (\ref{eq:new2}) implies that $ \uh(c)^{\oplus,\sharp}(\pi^\A)\le \Psi_\A(\uh(c)^\sharp(\pi^\A))<+\infty$  if  $\uh(c)^\sharp(\pi^\A)<+\infty.$}
\end{proof}

\begin{proof}[Proof of Theorem~\ref{theorem:newA}]

 The equivalence statement follows from  Lemma~\ref{lem:safe1bbbbbbaaaaa}.

``${\rm (L)}\Rightarrow{\rm (A4)}$'':  assumption~(L), stating boundedness below of the function $c,$  means that there exists a real number $\gamma\ge 0$ such that  $c(x,a)\ge -\gamma>-\infty$  for all $(a,b)\in\A\times\bb.$  This inequality implies that $c(a,b)+\gamma\ge c^+(a,b)$ for all $(a,b)\in\A\times\bb.$  By integrating  both sides of the last inequality in $\pi^\A\in \P(\A)$ and taking the supremum in $b\in\bb,$ we have $\uh(c)^\sharp(\pi^\A)+\gamma \ge \uh(c)^{\oplus,\sharp}(\pi^\A)$ for all $\pi^\A\in\P(\A).$
Therefore, if $\uh(c)^\sharp(\pi^\A)<+\infty,$ then $\uh(c)^{\oplus,\sharp}(\pi^\A)<+\infty,$ that is, assumption~(L) holds.

{The} remaining statements of Theorem~\ref{theorem:newA} follow from  Lemma~\ref{lem:safe1bbbbbbaaaaanew}.
\end{proof}


\begin{proof}[Proof of Theorem~\ref{th:exist_value}]
 {Condition~\eqref{ass:MAIN-NEW-A-WEAK} directly} implies formula \eqref{eq:valUvalS}. In view of Theorem~\ref{th:exvalPaper1}, the value of the game $v$ exists and the set $\P_v^\sharp(\A)$ is nonempty and convex. {Therefore, a}ccording to {Assumption~\ref{ass:MAIN-NEW-A} and} \eqref{eq:PvissafeEF}, $\P_v^\sharp(\A)\subset\P^S(\A).$  As follows from Definitions~\ref{defi:soleachequil} and \ref{Def5.8}, $\P_v^\sharp(\A)$ is the set of solutions for Player~I.
\end{proof}

The following corollary follows from Lemma~\ref{lem:safe1bbbbbbaaaaanew}.

\begin{corollary}\label{cor:safe1bbbbbbaaaaanew1}
 For a two-person zero-sum game
 $\{\A,\bb, c \},$
\begin{equation}\label{eq:new3B}
{\rm(L)}\Rightarrow{\rm(B1)}\Rightarrow{\rm(B2)}\Rightarrow{\rm(B3)}\Rightarrow{\rm(B4)}{\Rightarrow{\rm(B5)}}.
\end{equation}
\end{corollary}

\begin{proof}[Proofs of Theorem~\ref{theorem:newB} and Corollary~\ref{cor:safe1bbbbbbaaaaanew1}]
According to Remark~\ref{rem:new}, Theorem~\ref{theorem:newB} and Corollary~\ref{cor:safe1bbbbbbaaaaanew1}  follow respectively from Theorem~\ref{theorem:newA} and Lemma~\ref{lem:safe1bbbbbbaaaaanew} applied to the game $\{\bb,\A, -c^{\A\leftrightarrow\bb} \}$ introduced in Remark~\ref{rem:sim}.
\end{proof}

\subsection{Proofs of Theorem~\ref{teor:mainonestep}  and Corollary~\ref{cor:new}}\label{subsec:2}

\begin{proof}[Proof of Theorem~\ref{teor:mainonestep}]
 Theorem~\ref{th:exist_value} states the existence of the value. It also states that $\P_v^\sharp(\A)$ is the set of the solutions for Player I, $\P_v^\sharp(\A)\subset\P^S(\A),$ and $\P_v^\sharp(\A)$ is a nonempty convex compact subset of $\P(\A).$  Corollary~\ref{corth:exist_value} states that $\P_v^\flat(\bb)$ is the set of the solutions for Player II, $\P_v^\flat(\bb)\subset\P^S(\bb),$ and $\P_v^\flat(\bb)$ is a nonempty convex compact subset of $\P(\bb).$  As explained in the paragraph following formula \eqref{eq2.19EF}, $\P^\sharp_{{v}}(\A)\times \P^\flat_{{v}}(\bb)$
is the set of all solutions for the game.
\end{proof}

\begin{proof}[Proof of Corollary~~\ref{cor:new}]
The corollary directly follows from Theorems~\ref{theorem:newA}, \ref{theorem:newB}, and \ref{teor:mainonestep}. Indeed,
Theorems~\ref{theorem:newA} and \ref{theorem:newB} imply that Assumptions~{\ref{ass:MAIN-NEW-A} and \ref{ass:MAIN-NEW-B}} hold. Therefore, all assumptions of Theorem~\ref{teor:mainonestep}
hold because, in addition, the game $\{\A,\bb, c \}$ satisfy Assumptions~\ref{AsMain}(a1,a2,b1,b2).
\end{proof}

\subsection{Proofs of Propositions~\ref{propsec7A} and \ref{propsec7}}\label{subsec:3}
We start this subsection with some definitions and auxiliary lemmas.
We recall that for metric spaces $\X$ and $\Y$ a function $f:\X\times \Y\to \overline{\mathbb{R}}$ is called
\textit{$\K$-inf-compact on} $\X\times\Y,$
if for every  compact set $K\subset \X$ this function is inf-compact on $K\times \Y;$  see Feinberg et al. \cite[Definition 1.1]{Feinberg_et_al}, \cite[Definition~1]{OSCG}. A function $f:\X\times \Y\to \overline{\mathbb{R}}$ is called
 \textit{$\K$-sup-compact on} $\X\times\Y,$ if
 the function $-f$ is $\K$-inf-compact on $\X\times\Y;$ 
  see Feinberg et al. \cite[Definition~2]{OSCG}. 
 We would like to clarify that in this paper we consider $\K$-inf-compactness and $\K$-sup-compactness on the set $\X\times\Y,$ which the the graph
of the set-value mapping $\Phi:\X\to 2^\Y$ with $\Phi(x)=\Y$ for all $x\in\X,$ while in \cite{Feinberg_et_al} and \cite{OSCG} these definitions were considered for the graph of an arbitrary multifunction $\Phi$, where $\Phi:\X\to 2^\Y\setminus\{\emptyset\}$ in  \cite{Feinberg_et_al} and $\X:\Phi\to 2^\Y$ in \cite{OSCG}.  We observe that, if a function $f:\X\times \Y\to \overline{\mathbb{R}}$ is $\K$-inf-compact ($\K$-sup-compact) on $\X\times\Y,$ then for each $x\in\X$ the function $y\mapsto f(x,y)$ is inf-compact on $\Y.$ This follows from the observation that every singleton $K=\{x\},$  $x\in\X,$ is compact.

\begin{lemma}{\rm(Feinberg et al. \cite[Lemma~2]{OSCG})}\label{k-inf-compact}
The function $f:\X\times \Y\to \overline{\mathbb{R}}$ is
$\K$-inf-compact on $\X\times \Y$ if and only if the following two assumptions hold:
\begin{itemize}
\item[{\rm
(i)}] $f:\X\times \Y\to \overline{\mathbb{R}}$ is lower semi-continuous;
\item[{\rm(ii)}] if a sequence $\{x^{(n)} \}_{n=1,2,\ldots}$ with values in $\X$
converges in $\X$ and its limit $x$ belongs to $\X,$ then each sequence $\{y^{(n)}
\}_{n=1,2,\ldots}\subset\Y$ satisfying
the condition that the sequence $\{f(x^{(n)},y^{(n)}) \}_{n=1,2,\ldots}$ is
bounded above, has a limit point $y\in \Y.$
\end{itemize}
\end{lemma}

The following lemma is the main technical fact in this subsection.

\begin{lemma}\label{lem:exa1}  Let
$\A=\bb:=\mathbb{R}_+$ and $c(a,b):=\varphi(a-b)$ for each $a,b\in \mathbb{R}_+,$ where $\varphi:\mathbb{R}\to\mathbb{R}$ is a continuous function.
Then the following statements hold:
\begin{itemize}
\item[(i)] if $\varphi(s)\to+\infty$ as $s\to+\infty,$ then the function $(b,a)\mapsto c(a,b)$ is $\K$-inf-compact on $\bb\times\A;$
\item[(ii)] if $\varphi(s)\to-\infty$ as $s\to-\infty,$ then the function $(a,b)\mapsto c(a,b)$ is $\K$-sup-compact on $\A\times\bb;$
\item[(iii)] if $\varphi(s)=\varphi_1(s)+\varphi_2(s)$ for each $s\in\mathbb{R},$ where $\varphi_1:\mathbb{R}\to\mathbb{R}$ is increasing and $\varphi_2:\mathbb{R}\to\mathbb{R}$ is bounded, then  assumptions~(A1) and (B1) from Theorems~\ref{theorem:newA} and \ref{theorem:newB} hold;
\item[(iv)] if there exist $s_*<0<s^*$ such that $\varphi(s_*)>\varphi(s^*),$ then the game $\{\A,\bb,c\}$ has  no pure solution.
\end{itemize}
\end{lemma}

\begin{proof}
(i) We verify the conditions of Lemma~\ref{k-inf-compact} to prove $\K$-inf-compactness of the function $c(a,b)=\varphi(a-b).$ This function is continuous, and therefore it is lower semi-continuous.   Consider a sequence $\{b^{(n)}\}_{n\ge 1}$ that converges to $b\in \bb$ and a sequence $\{a^{(n)}\}_{n\ge 1}\subset \A$ such that  $\{\varphi(a^{(n)}-b^{(n)})\}_{n\ge 1}$ is bounded above. Since the sequence $\{b^{(n)}\}_{n\ge 1}\subset \mathbb{R}_+$ converges, it is bounded.   Since the sequence  $\{\varphi(a^{(n)}-b^{(n)})\}_{n\ge 1}$ is bounded above, then the continuity of the function $\varphi:\mathbb{R}\to\mathbb{R}$ on $\mathbb{R}$ and the property $\varphi(s)\to+\infty$ as $s\to+\infty$ imply that the sequence $\{a^{(n)}-b^{(n)}\}_{n\ge 1}$ is bounded above.  Thus, the sequence $\{a^{(n)}\}_{n\ge 1}\subset \mathbb{R}_+$ is bounded above and therefore it is bounded.  Therefore, the sequence $\{a^{(n)}\}_{n\ge 1}$ has an accumulation point $a\in \A$. Thus, the assumptions of Lemma~\ref{k-inf-compact} are verified, and the function $c$ is $\K$-inf-compact.

(ii)  This statement follows from (i) applied to the game $\{\bb,\A,-c^{\A\leftrightarrow\bb}\}.$

(iii) First, we prove that assumption~(B1) holds. Let the function $\varphi$ be the sum of the functions $\varphi_1$ and $\varphi_2$ described in the statement. Then for each $b\ge0$
\[
c^\flat(b)=\inf_{a\ge0}\{\varphi_1(a-b)+\varphi_2(a-b)\}\ge \inf_{a\ge0}\varphi_1(a-b)+\inf_{a\ge0}\varphi_2(a-b)
\]
\[
=\varphi_1(-b)+\inf_{a\ge0}\varphi_2(a-b)=c(0,b)+\inf_{a\ge0}\varphi_2(a-b)\ge c(0,b)-B,
\]
where $B>0$ is a constant such that $|\varphi_2(s)|\le B$ for each $s\in\mathbb{R}.$ We note that the second equality holds because the function $\varphi_1$ is increasing. Therefore, for each $b\ge0$
\[
c^-(0,b)\le \frac12 c^-(0,b)\le \frac12c^\flat(b)+\frac{B}2.
\]
This implies that assumption~(B1) holds.

Second, assumption~(A1) holds because it is equivalent to assumption~(B1) for the game $\{\bb,\A, -c^{\A\leftrightarrow\bb} \},$ which holds because the real function $\varphi_{1}$ is increasing if and only if
the real function $s\mapsto -\varphi_{1}(-s)$ is increasing, and the function $\varphi_{2}$ is bounded if and only if
the function $s\mapsto -\varphi_{2}(-s)$ is bounded.

(iv) There exist $s_*,s^*\in\mathbb{R}$ such that $s_*<0<s^*$ and $\varphi(s_*)>\varphi(s^*).$ Then for each $a,b\ge0$
\[
c^\flat(b)=\inf_{a^*\ge 0} \varphi(a^*-b)\le \varphi(s^*),\quad c^\sharp(a)=\sup_{b^*\ge 0}\varphi(a-b^*)\ge \varphi(s_*).
\]
Therefore,
\[
\sup_{b\ge 0}c^\flat(b)\le \varphi(s^*)< \varphi(s_*)\le \inf_{a\ge 0}c^\sharp(a),
\]
that is, the game $\{\A,\bb,c\}$ has no  pure solution. 
\end{proof}

\begin{proof}[Proof of Proposition~\ref{propsec7A}.] (a) In view of Lemma~\ref{lem:exa1}(i), the function $(b,a)\mapsto c(a,b)$ is $\K$-inf-compact on $\bb\times\A.$  This implies Assumptions~\ref{AsMain}(a1,a2).
Statement (b) follows from Lemma~\ref{lem:exa1}(i,iii). Statement (c) follows from Lemma~\ref{lem:exa1}(i-iii).
\end{proof}

\begin{proof}[Proof of Proposition~\ref{propsec7}](a) This statement follows from Proposition~\ref{propsec7A}(a) and Lemma~\ref{lem:exa1}(iv).

(b) Let us consider three cases (c1--c3).

(c1) Let $M$ be even and $\alpha_M>0.$ Then condition~(v) from Definition~\ref{defi:game}
does not hold because the function $b\mapsto c(a,b)$ is not bounded above on $\mathbb{R}$ for each $a\ge0.$
 Since the function $\varphi$ is
bounded  below on $\mathbb{R},$ the value $\uh(c)(\pi^\A,\pi^\bb)$ is well-defined for  all $(\pi^\A,\pi^\bb) \in\P(\A)\times \P(\bb)$ and
\begin{equation}\label{eq:newformula}
\sup_{\pi^\bb\in\P(\bb)}\inf_{\pi^\A\in\P(\A)}\uh(c)(\pi^\A,\pi^\bb)=\inf_{\pi^\A\in\P(\A)}\sup_{\pi^\bb\in\P(\bb)} \uh(c)(\pi^\A,\pi^\bb)=+\infty.
\end{equation}
 Indeed, if we set $\pi^\bb(B):=\frac{2}{\pi}\int_{B}\frac{1}{1+b^2}db$ for each $B\in\B(\bb),$ then for all $a\in\A$
\[
\hat{c}(a,\pi^\bb)=\frac{2}{\pi}\int_{\mathbb{R_+}}\frac{\varphi(a-b)}{1+b^2}db=+\infty.
\]
 Therefore,
$\hat{c}^\flat(\pi^\bb)=\inf_{\pi^\A \in \P(\A)}\hat{c}(\pi^\A,\pi^\bb)=\inf_{a\ge0}\hat{c}(a,\pi^\bb)= +\infty,$ and
\[
+\infty\le \hat{c}^\flat(\pi^\bb)\le
\sup_{\pi^\bb_*\in\P(\bb)}\hat{c}^\flat(\pi_*^\bb)=\sup_{\pi^\bb\in\P(\bb)}\inf_{\pi^\A\in\P(\A)}\uh(c)(\pi^\A,\pi^\bb)\le\inf_{\pi^\A\in\P(\A)}\sup_{\pi^\bb\in\P(\bb)} \uh(c)(\pi^\A,\pi^\bb).
\]
 Thus equalities (\ref{eq:newformula}) hold.

(c2) If $M$ is even and $\alpha_M<0,$ then
condition~(iv) from Definition~\ref{defi:game} does not hold because the function $a\mapsto c(a,b)$ is not bounded below on $\mathbb{R}$ for each $b\ge0.$ Since the function $\varphi$ is
bounded above on $\mathbb{R},$  the value $\uh(c)(\pi^\A,\pi^\bb)$ is well-defined for  all $(\pi^\A,\pi^\bb) \in\P(\A)\times \P(\bb).$ Moreover, by the symmetric reasonings, which follow from case (c1),
\[
\sup_{\pi^\bb\in\P(\bb)}\inf_{\pi^\A\in\P(\A)}\uh(c)(\pi^\A,\pi^\bb)=\inf_{\pi^\A\in\P(\A)}\sup_{\pi^\bb\in\P(\bb)} \uh(c)(\pi^\A,\pi^\bb)=-\infty.
\]

(c3) If $M$ is odd and $\alpha_M<0,$ then
conditions~(iv,v) from Definition~\ref{defi:game} do not hold. Moreover,  the lower value for this game in pure strategies equals  $-\infty,$ and the upper value for this game in pure strategies equals  $+\infty.$
\end{proof}

\medskip

{\bf Acknowledgements.}  Research of the first author was partially supported by NSF Grant CMMI-1636193.  The authors thank William D. Sudderth for his valuable comments on von Neumann's
and Sion's minimax theorems.

\end{document}